



\font\bm=msym10

\def\m#1{\hbox{\bm #1}}




\vsize=24.5truecm
\hsize=6.5truein


\nopagenumbers
\headline={\ifnum\pageno=1\hfil\else\hfil\rm\folio\fi}

\def\t#1{\widetilde{#1}}
\def\h#1{\widehat{#1}}
\def\b#1{\overline{#1}}
\def\oa#1{\buildrel #1\over\rightarrow}
\def\blowup#1{\,\t#1\oa\pi #1\,}
\def\isoto{\rightarrow \kern-13pt \lower.83pt\hbox{${}^\sim$}~\,}
\def\onto{\rightarrow\kern-8.8pt\rightarrow}
\def\norm#1{|\kern-1.72pt{|}\,#1\,|\kern-1.72pt{|}}
\def\sqr#1#2{{\vbox{\hrule height.#2pt \hbox{\vrule width.#2pt
  height#1pt \kern#1pt \vrule width.#2pt}\hrule height.#2pt}}}
\def\qed{\qquad \sqr 7 3}

\def\d{\partial}
\def\db{\bar\partial}
\def\dbd{\bar\partial\partial}
\def\ddb{\partial\bar\partial}


\centerline{\bf Blowups and gauge fields}

\bigskip
\centerline{N.~P.~Buchdahl \footnote*{Research partially
supported by NSF Grant \#8900878 \hfill\break A.M.S.  (1980)
subject classification ~~32L10,14F05, 53B35}}
\centerline{Department of Pure Mathematics}
\centerline{University of Adelaide}
\centerline{Adelaide, Australia 5005}
\centerline{e-mail:
nbuchdah@maths.adelaide.edu.au}
\vskip.75in

\noindent{\bf 0.\quad Introduction.}

\bigskip
The purpose of this paper is to investigate the
relationship between between stable holomorphic vector bundles on
a compact complex surface and the same such objects on a
modification (blowup) of the surface.  In large part, the paper
is a continuation of the work in [B2] where it was shown that a
holomorphic bundle on a compact complex surface admits an
irreducible Hermitian-Einstein connection if and only if the
bundle is stable, the notions of {\it stability\/} and {\it
Hermitian-Einstein\/} both being with respect the same
$\dbd$-closed positive $(1,1)$-form: this is a generalization of
Donaldson's result [D2] where it is assumed that the form is
$d$-closed and defines an integral cohomology class (i.e., the
algebraic case).

Much of the underlying motivation for this work comes from its
potential applications to topology (though no such applications
are considered here). Donaldson has proved fundamental results on
the  topology of smooth 4-manifolds by defining topological
invariants  of moduli spaces of solutions of the anti-self-dual
Yang-Mills  equations on a Riemannian 4-manifold, and showing
that these are  differential invariants of the manifold itself.
In the case of an  algebraic surface with Hodge metric, the
result of [D2] mentioned  above identifies the Yang-Mills moduli
spaces with moduli spaces of  stable vector bundles, and these
spaces and/or invariants can be computed using techniques
standard in complex analysis.  In this  way Donaldson has been
able to prove some of his most remarkable  results [D4], [D5],
and others have built upon his work ([FM1],  [FM2], [K1], [OV],
to cite just a few. See also [FM3] for a comprehensive account of
developments, and [FS] for a calculation of the Donaldson
invariants of a blownup $4$-manifold in terms of those of the
manifold).  This interaction between real analysis, complex
analysis and topology provides a rich area for investigation, and
parts of this paper are directly concerned particularly with the
interplay between the real and the complex analysis.

The proof of the main result of [B2] is a modification of that
given by Donaldson [D1] to prove the same theorem in the case of
Riemann surfaces.  The differences in the proofs arise from the
appearance of certain singularities in the two-dimensional case,
and a successful way around these singularities is to blow up the
surface and pull back.  In so doing, various relationships
between bundles and sheaves on the surface and its blowups are
uncovered, and these relationships turn out to be directly
related to other aspects of gauge theory and/or complex analysis
which are themselves of independent interest.

A study of degenerating sequences of stable bundles on the
projective plane leads naturally to conjecture whether blowups
can be used to compactify moduli spaces of stable bundles in
general.  In [B4], it is shown that sequences of stable bundles,
identified with sequences of Hermitian-Einstein connections have
convergent subsequences after pulling back to blowups, at least
when weak limits are stable. This leads to a the definition of a
natural topology on moduli spaces stable bundles over a surface
and its blowups, and the proof of the compactness of the generic
such space is presented here.

\medskip
The paper is organized as follows: \S 1 introduces
notation, definitions and central background material, and gives
some useful lemmas concerning ``invariants" of stable holomorphic
bundles.  In \S 2 a local description and characterization of
bundles on the blowup of the ball in $\m C^2$ at the origin is
given.   A holomorphic version of Taubes' ``cut-and-paste"
construction for gauge fields [T] is given, enabling a global
description of bundles on the blowup of an arbitrary complex
surface in terms of bundles on the original surface.  Also
included in this section is a short discussion of the
relationship to---and between---associated constructions of Serre
and Schwarzenberger.

Questions of stability are considered in the third section from a
purely complex-analytic viewpoint, and a detailed description of
the conditions required for bundles on a blowup $\blowup X$ to be
stable is given.  Of course, the metric on $\t X$ with respect to
which stability is considered must be specified, and the
fore-mentioned conditions are as much on the metric as on the
bundles themselves.  From a real analytical view-point, it turns
out the correct metrics are those formed from the original metric
on $X$ together with the Fubini-Study metric on $\b{\m {P}}_2$
combined so as to ``stretch-out" the neck of the connected sum
$\t X \simeq_{diffeo} X\#\b{\m {P}}_2$, and the main result of
this section shows that once the neck is sufficiently long, the
moduli spaces effectively become independent of the metric.

The analysis in the third section encounters pathological sheaves
which are semi-stable but not stable.  Using the cut-and-paste
method, a mechanism for ``stabilising" such sheaves is given in
\S4.  A similar method also provides a simple way to
desingularise singular points in moduli spaces.

For a bundle in one of the ``stable" moduli spaces of \S3, the
behaviour of the corresponding sequence of Hermitian-Einstein
connections as the metrics degenerate is investigated in \S5.

In the last section, the issue of compactness for moduli spaces
of stable bundles is considered.  In [B4] it is shown that after
sufficiently many blowups and pull-backs, sequences of stable
bundles of bounded topology and degree have strongly convergent
subsequences, where stable bundles are identified with
irreducible Hermitian-Einstein connections.   A natural candidate
for a compactification of a moduli space of stable bundles as
presented in [B4] is shown, under generic conditions in the
arbitrary rank case, and in general for the rank 2 case, to be a
compact space;  some other simple properties of this space are
also considered.

\bigskip
\bigskip
\noindent{\bf Acknowledgments:} \quad Parts of
this paper were written in 1988 when I  was a visiting member of
the Institute for Advanced Study,  Princeton, N.~J. supported by
the National Science Foundation. Other parts were written  while
a member of  the Department of Mathematics at Tulane University,
New Orleans  with the support of NSF Grant \#DMS 8900878. The
work was completed while visiting the Department of Mathematics
at the University of Nantes, to which institution I am very
grateful for its hospitality.  I am also grateful to the \'Ecole
Polytechnique, Paris for its hospitality during the same period,
where further additions to the paper were made.

\bigskip
\bigskip
\bigskip

\noindent {\bf 1.  \quad Preliminaries.}

\bigskip
\medskip
The purpose of this section is to re-cap on,
and to expand upon the basic notation, definitions and results of
[B2].  Further details can be found in that reference.
\bigskip
Let $X$ be a compact complex surface and let $\omega$ be a $\bar
\partial\partial$-closed positive (1,1)-form on $X$: it is a
theorem of Gauduchon [Gau] that every positive $(1,1)$-form has a
unique positive conformal rescaling such the rescaled form is
$\dbd$-closed and gives the same volume $V := Vol(X,\omega) :=
{1\over 2}\int_X \omega^2$.  With such a form $\omega$, the
degree $deg(L) = deg(L,\omega)$ of a holomorphic line bundle $L$
on $X$ is unambiguously defined by the formula
$$
deg(L) \; := \;
{i\over 2\pi} \int_X f_L \wedge \omega\;,
$$
where $f_L$ is the
curvature of any hermitian connection on $L$.  The degree depends
only on $c_1(L)$ if and only if $b_1(X)$ is even, and when this
is the case $\omega$ is cohomologous modulo the image of
$\partial + \bar\partial$ to a closed form which itself is unique
up to the image of $\dbd$; ([B2], Proposition~2).

If $E$ is a holomorphic $r$-bundle on $X$, set $deg(E) :=
deg(det\,E)$ and $\mu (E) := deg(E)/r$; the latter is called the
normalized degree  or {\it slope\/} of $E$.  A hermitian
connection on $E$ is {\it  Hermitian-Einstein\/} if the curvature
$F$ satisfies $\h F =  i\lambda \bf 1$ where $\h F := *\,(\omega
\wedge F) =: \Lambda F$,  $\lambda = (-2\pi / V )\!\cdot\! \mu
(E)$ and $\bf 1$ is the  identity endomorphism of $E$.  The
bundle $E$ is ({\it semi-}) {\it  stable\/} if $\mu ({\cal S}) <
(\le)~\mu (E)$ for every coherent  subsheaf ${\cal S} \subset E$
with $0 < rank({\cal S}) < r$.  As  mentioned in the
introduction, the main result of [B2] is that a  bundle admits an
irreducible Hermitian-Einstein connection if and  only if it is
stable, this generalizing the same result proved by  Donaldson
[D2] in the case that $(X,\omega)$ is algebraic.  A  bundle
admitting a Hermitian-Einstein connection is a direct sum of
stable bundles all of the same normalized degree; i.e., is {\it
quasi-stable\/}.

If $E$ has a Hermitian-Einstein connection with curvature $F$,
the equation $\omega \wedge \big( F - {1\over r}tr F\, {\bf 1}
\big) = 0$ and the skew-Hermitian property of $F$ give $tr\big( F
- {1\over r}tr F \,{\bf 1} \big)^2 = | F - {1\over r}tr F \, {\bf
1}|^2\,dV$.  Since the former $4$-form is a representative for
the characteristic class $8\pi^2\big( c_2 - {r-1\over
2r}c_1^2\big)(E)$, this motivates defining the {\it charge\/} of
$E $, $C(E)$, for an arbitrary $r$-bundle $E$ by the formula
$$
C(E) := \big( c_2 - {r-1\over 2r}c_1^2\big)(E) \;=\; {1\over
8\pi^2}\int_X tr\big( F - {1\over r}tr\, F \,{\bf 1} \big)^2
\quad .\eqno (1.1)
$$
This number is non-negative for any bundle
admitting a Hermitian-Einstein connection, and when this is the
case, is identically zero only if the induced Hermitian-Einstein
connection on the adjoint bundle is flat; (cf.~[L]).  Note that
the charge is invariant under tensoring by line bundles:
$C(E\otimes L) = C(E)$ for any such $ L $.  In general,
$C(E\otimes A) = aC(E) + rC(A)$, where $a,r$ are the ranks of
$A,E$ respectively.
\medskip
Recall that a coherent analytic
sheaf ${\cal S}$ is torsion-free if and only if the canonical
morphism ${\cal S} \to {\cal S}^{**}$ is injective, and ${\cal
S}$ is by definition reflexive if this map is an isomorphism;
recall also that the singularity sets of such sheaves are of
codimension at least 2 and 3 respectively; ([OSS], II.1.1).  For
exact sequences $\;0 \to {\cal A} \to {\cal B} \to {\cal C} \to 0
\;$ of locally free sheaves on $X$ it is easy to check that the
charges are related by
$$
C({\cal B}) = C({\cal A}) + C({\cal C})
- {b\over 2ac}\big[ {a\over b}c_1({\cal B}) - c_1({\cal
A})\big]^2\; ,\eqno(1.2)
$$
where $a$, $b$ and $c$ are the ranks
of ${\cal A}$, ${\cal B}$ and ${\cal C}$ respectively.  The
definition of charge extends to torsion-free sheaves ${\cal S}$
of rank $r$ by means of the formula
$$
C({\cal S}) := C({\cal
S}^{**}) + h^0({\cal S}^{**}/{\cal S}) \;, \eqno(1.3)
$$
which is
consistent with a definition of $c_2({\cal S})$ extending that of
the Chern character on bundles in such a way that the
Hirzebruch-Riemann-Roch formula
$$
h^0({\cal S})-h^1({\cal
S})+h^2({\cal S}) = \chi({\cal S})=-C({\cal S})+{1\over
2r}c_1^2({\cal S}) +{1\over 2}c_1({\cal S})\cdot c_1(X) + r\chi
({\cal O}_X)
$$
remains valid.  If ${\cal C}$ is  only
torsion-free, it follows from this definition that (1.2) remains
valid and this in turn implies that the formula (1.2) holds for
arbitrary torsion-free sheaves ${\cal A}$, ${\cal B}$ and ${\cal
C}$. Note that a torsion-free sheaf is (semi-)stable iff its
double-dual is.

\medskip
If $b_1(X)$ is odd, the intersection form on
$H^2(X,\m{R})$ restricted to $H^{1,1}(X)$ is negative definite
([BPV], Theorem~ IV.2.13) and the last term on the right in (1.2)
therefore contributes positively to the sum.  If $b_1(X)$ is
even, the intersection form restricted to $H^{1,1}(X)$ has one
positive eigenvalue and the rest are all negative.  In either
case, $\omega$ defines a positive definite hermitian form on
$H^{1,1}(X)$ by setting $|\!| f |\!|^2 := V^{-1}|(f,\omega)|^2 -
(f,f) $, where $(f,g) := \int_X \bar f \wedge g $; (recall $V =
(\omega,\omega)/2$ throughout).  Equation (1.2) can therefore be
written
$$
C({\cal B}) \; = \;C({\cal A}) + C({\cal C}) + {b\over
2ac}|\!| {a\over b}c_1({\cal B})- c_1({\cal A})|\!|^2 - {b\over
2ac} {\nu_ {\cal B} ({\cal A})^2\over V}\; ,\eqno (1.4)
$$
where
$\nu_ {\cal B} ({\cal A}) := a\big[ \mu({\cal B}) - \mu ({\cal
A}) \big] $.  By induction on rank, it follows the charge is
non-negative for any torsion-free semi-stable sheaf.  Note that
if $b_1(X)$ is odd it follows by induction from (1.2) (and the
existence of Hermitian-Einstein connections on stable bundles)
that the charge  is non-negative for {\it any\/} torsion-free
coherent analytic sheaf, semi-stable or otherwise.

The function $\nu_{\bullet}(*)$ plays an important role in the
proof of the main result of [B2].  It has a number of simple but
useful properties, three of which are summarised for convenience
in the following lemma.

\bigskip
\noindent {\bf Lemma~1.5.} {\sl
\medskip
\item{(a)~} If
$$
\matrix{0&\longrightarrow&{\cal A}&\longrightarrow&{\cal
B}&\longrightarrow &{\cal C}&\longrightarrow&0\cr
&&\uparrow&&\uparrow \vrule height15pt width0pt&&\uparrow&&\cr
0&\longrightarrow&{\cal A}'&\longrightarrow&{\cal B}'\vrule
height15pt width0pt&\longrightarrow &{\cal
C}'&\longrightarrow&0\cr\cr}
$$
is a commutative diagram with
exact rows such that the vertical arrows are inclusions, then
$$
\nu_{\cal B} ({\cal B}') = \nu_{\cal A} ({\cal A}') + \nu_{\cal
C} ({\cal C}') + \big( {a'\over a} - {c'\over c}\big) \nu_{\cal
B} ({\cal A})\;; \eqno (1.6)
$$
\item{(b)~} If ${\cal A}$ and
${\cal B}$ are locally free and if ${\cal B}$ is stable and
${\cal A} \subset {\cal B}$ minimizes $\nu _{\cal B}$ over all
proper non-zero subsheaves of ${\cal B}$ then ${\cal A}$ is
stable and in addition, the quotient ${\cal C} = {\cal B}/{\cal
A}$ is both torsion-free and stable;
\item{(c)~} If $E$ is a
holomorphic bundle equipped with a Hermitian metric and $A
\subset E$ has torsion-free quotient $C$ then off the singular
set of $C$ the second fundamental form $\beta \in
\Lambda^{0,1}\otimes Hom(C,A)$ of the induced Hermitian
connection lies in $L^{2p}(X)$ for any $p < 2, $ and
$$
\nu_E(A)
= -{1\over 2\pi} \int_X tr_A\big(i\,F(E)+ \lambda_E\, {\bf
1}\,\omega \big)\wedge\omega \; +\; |\!|\beta|\!|^2_{L^2(\omega)}
\;.  \eqno (1.7)
$$
\qed}

\bigskip
\noindent The proof of (c) is given in [B2] (remark (f)
p.634).  Part (b) is the same as Lemma~2 of the same reference;
the proof follows immediately from (1.6) which itself is a
straight-forward calculation.  The existence of such ${\cal A}
\subset {\cal B}$ minimizing $\nu_{\cal B}$ when the latter is
stable is proved in Lemma~4 of [B2], which provides one of the
key steps in the proof of the main result there by enabling the
argument to proceed by induction on rank; (it is also proved in
that lemma that there always exists ${\cal A} \subset {\cal B}$
maximizing $\mu$ over admissible subsheaves regardless of the
stability or otherwise of ${\cal B}$).  For a stable bundle $E$
it follows that $\nu_E(*)$ is bounded above and away from $0$ on
the set of subsheaves of $E$, from which it follows immediately
that stability is an open condition on the metric.

When $b_1(X)$ is even, this bound on the slopes of subsheaves can
be made uniform in $E$ as the  next result shows:

\bigskip
\noindent {\bf Lemma~1.8.} \sl Suppose that $ b_1(X)$ is
even.  For any $r_0,C_0 > 0$ there exists $\delta_0 =
\delta_0(r_0,C_0) > 0$ with the following property: if $E$ is a
semi-stable torsion-free sheaf of rank $r \le r_0$ and charge
$C(E) \le C_0$ admitting a subsheaf $A \subset E$ with $\nu_E(A)
< \delta_0 $, then $\nu_E(A) = 0$.  \rm
\medskip
\noindent {\bf
Proof: } Since $C(E^{**}) \le C(E)$ and $\nu_{E^{**}}(A^{**}) =
\nu_E(A)$ it suffices to prove the result with ``torsion-free" in
the hypotheses replaced by ``locally free".  For subsheaves $A
\subset E$ of the same rank as $E$ the result follows from
Corollary~2 of [B2], so it can also be supposed that all such
subsheaves have rank strictly less than that of the ambient
bundle.

If there is no stable bundle of rank $ \le r_0$ and charge $\le
C_0$ which admits a proper non-zero reflexive subsheaf, then
$\nu_E$ is identically $0$ for all such bundles.  Otherwise,
there is a sequence $\{E_i\}$ of bundles admitting such
subsheaves $A_i \subset E_i$ with $\{\nu_{E_i}(A_i)\}$ strictly
decreasing.  By Lemma~4 of [B2], it can be assumed that each
$A_i$ minimizes $\nu_{E_i}$ over the proper non-zero subsheaves
of $E_i$, and by Lemma~1.5(b), both $A_i$ and the quotient $C_i
:= E_i/A_i$ are torsion-free and stable.  Since $ C(A_i)$ and $
C(C_i)$ are therefore non-negative and $\nu_{E_i}(A_i) =
|\nu_{E_i}(A_i)|$ is decreasing, (1.4) and the bound on the rank
and charge of $E_i $ give a uniform bound on $\norm{a_i c_1(E_i)
- r_ic_1(A_i) }$, so there is a subsequence with $ r_i $ and $a_i
c_1(E_i) - r_ic_1(A_i)$ constant.  Since $b_1(X) $ is even, the
degree is topological so $\nu_{E_i}(A_i)$ is constant on the
subsequence, implying that the original sequence is finite.  It
follows that there exists $\delta > 0$ such that $\nu_E(A) \ge
\delta $ for any proper non-zero subsheaf $A$ of a stable
torsion-free sheaf $E$ of rank $\le r_0$ and charge $ \le C_0$.
Set $\delta_0 := \delta $.

If now $E$ is semi-stable (of rank $\le r_0$ and charge $ \le
C_0$) and $D \subset E$ satisfies $\nu_E(D) < \delta_0$, then $E$
cannot be stable so there exists non-trivial $A \subset E$ with
torsion-free quotient $C$ such that $\mu (A) = \mu (E)$.  Since
$E$ is semi-stable, so too are both $A$ and $C$ and moreover they
also satisfy the hypotheses of the lemma by (1.4).  If $ C'$ is
the image of the composition $D \to E \to C$ and $A'$ is the
kernel, the result follows from (1.6) using induction on
$rank(E)$.  \qed
\bigskip
\noindent {\bf Remark: } The result is
false if $b_1(X)$ is not even: moduli spaces of stable
$2$-bundles with trivial determinant on a Hopf surface are
explicitly computed in [BH], and the description there shows that
there are stable $2$-bundles of charge $1$ possessing subsheaves
of degree arbitrarily close to $0$.

\bigskip
\bigskip
\bigskip

\noindent{\bf 2.\quad Vector bundles on a blowup.}

\bigskip
The purpose of this section is to investigate the nature
of holomorphic vector bundles on a neighbourhood of a blown-up
point in a complex surface.  Two approaches are taken:  the first
starts from the splitting type of such bundles on the exceptional
divisor, whereas the second is more global in nature, identifying
such bundles with a class of bundles on the complex projective
plane.
\medskip
Let $Y$ be a discrete set of points in a complex
surface $X$ and let $\blowup X$ be the blowup of $X$  along $Y$.
The exceptional divisor
$\t Y = \pi^{-1}(Y) \subset \t X$ is defined by a section of a
certain holomorphic line bundle, and since this line bundle
restricts to ${\cal O}(-1)$ on each component of $\t Y$ the
notation ${\cal O}(-1)$ will  be used to denote (the sheaf of
sections of) this line bundle. If  ${\cal I}_Y \subset {\cal
O}_X$ denotes the ideal sheaf of $Y$ and ${\cal N}_Y=({\cal
I}_Y/{\cal I}_Y^2)^*$ is the normal bundle of $Y$ in $X$, then
it is straightforward to show that the direct images of the
sheaves ${\cal O}(n)$ under $\pi$ are canonically given by
$$
\eqalign{\pi_*{\cal O}_{\t X}(n) =&  \cases{{\cal I}_Y^n&\quad if
$n\ge 0$ \cr
{\cal O} &\quad if $ n\le 0$} \cr\cr
\pi_*^1{\cal
O}_{\t X}(n) =&   \cases{$0$ &\quad if  $n \ge -1$ \cr
{\rm
det\;}{\cal N}_Y &\quad  if $n = -2$ \cr
{\cal E}xt^1_{\cal
O}({\cal I}_Y^{-n-1},{\cal O}) & \quad if $ n \le -2$}}\eqno(2.1)
$$
(where $\pi_*^q = R^q\pi_*$ denotes the $q$-th direct image
under $\pi$) with all other direct images  vanishing.

\medskip
Let $ U \subset X$ be a small ball in $X$ such that
$U\cap Y$ is the singleton $x_0$,  let $\blowup U$ be the blowup
of $U$ at $x_0$, and let $L_0= \pi^{-1}(x_0)$ be the exceptional
line.  If $\t E$ is a holomorphic $r$-bundle on $\t U$, then the
restriction of $\t E$ to $L_0$ splits as a sum of line bundles,
and the nature of this splitting determines much about the bundle
itself, as will be demonstrated in the results which follow.

\bigskip
\noindent{\bf Lemma~2.3.}\quad {\sl Suppose $ \t
E\!\!\mid _{L_0} \, = \oplus^r_{i=1} {\cal O}(a_i)$.
\medskip
\itemitem{\rm (a)\quad} If $a_i \le 0 $ for all $i$ then $\pi_*\t
E$ is locally free;
\itemitem{\rm (b)\quad} If $\pi_*\t E$ is
locally free then $\sum_i a_i \le 0$;
\itemitem{\rm (c)\quad} If
$a_i=0$ for all $i$ then $\t E = \pi^*\pi_*\t E $ and is trivial
on $\t U$;
\itemitem{\rm (d)\quad}  If $a_i > -2$ for all $i$
then $\pi_*^1\t E = 0$.}

\medskip
\noindent{\bf Proof: } The space $\t U$ can be viewed as
a closed  subspace of $U \times \m{P}_1$, defined by a section of
the line bundle ${\cal  O}(1)$. The bundle $\t E$ on $\t U$ can
be extended (non-uniquely) to a bundle $\t E'$ on $U\times \m
P_1$ simply by extending a transition function on the
intersection of a pair of Stein sets covering $\t U$. Thus there
is an exact sequence
$$
0 \to \t E'(-1) \to \t E' \to \t E\to 0
\eqno (2.4)
$$
on $U\times \m P_1$, (where the notation $\t E(n)$
denotes $\t E\otimes {\cal O}(n)$ throughout).

If $a_i< 0$ for all $i$ then $H^0(\{x\}\times \m P_1,\t E')=0$
for $x=x_0$ and hence for  all $x$ in a neighbourhood of $x_0$ by
semi-continuity of cohomology.  From the base-change theorem
([BS], Theorem~3.4)  it follows that  $\pi_*\t E'= 0= \pi_*\t
E'(-1)$  and the sheaves $\pi^1_*\t E'$ and $ \pi^1_*\t E'(-1)$
are locally free on $U$. From the direct image of (2.4)
$$
0\to
\pi_*\t E \to \pi^1_*\t E'(-1) \to \pi^1_*\t E' \to \pi^1_*\t E
\to 0
$$
it therefore follows from  Lemma~II.1.1.10 of [OSS] that
$\pi_*\t E$ is reflexive, and hence locally free. Moreover,
taking direct images of the exact sequence $\t E\otimes (0 \to
{\cal O}_{\t U}(1) \to {\cal O}_{\t U} \to {\cal O}_L \to 0) $ on
$\t U$ and using the fact that $\pi_*(\t E\!\!\mid_{L_0})$
vanishes, it follows that $\pi_*\t E(1) = \pi_*\t E$  is also
reflexive, which proves (a).

To prove (b), suppose $\det \t E = {\cal O}(a)$ for some $a>0$.
If $\pi_*\t E$ is locally free, then by (2.1) so too is $\pi_*\t
E''$  for $\t E'' := \t E\oplus{\cal O}(-a)$, and  therefore
$\pi^*\pi_*\t E''$ is a bundle on the blowup.  The canonical
sheaf homomorphism $\pi^*\pi_*\t E''\to \t E''$ is an isomorphism
off the exceptional divisor, but since both bundles have trivial
determinant the homomorphism must in fact be an isomorphism
everywhere.  This is not possible since the pull-back is trivial
on  the exceptional divisor, whereas $\t E''$ is not.

If $\t E$ is trivial on $L_0$, then $\pi_*\t E$ is locally free
by (a) and the argument of the last paragraph shows that $\t E
\simeq \pi^*\pi_*\t E$, proving (c).

Finally, if $a_i > -2 $ for all $i$ then $H^1(\{x\}\times \m
P_1,\t E')$ vanishes at $x=x_0$, and therefore for all $x$ in a
neighbourhood of $x_0$, again by semi-continuity of cohomology.
Taking direct images of (2.4)  this time shows that $\pi^1_*\t E$
vanishes near $x_0$.  \qed
\bigskip

With $\t E$ as in Lemma~2.3, suppose that $a_1 \le \dots \le
a_r$. The obstruction to extending a section in $\Gamma(L_0,\t
E(-a_1))$ to $\t U$ lies in  $H^1(\t U,\t E(1-a_1))$, a group
which vanishes by part (d) of the lemma. More generally, any
finite number of such sections which are independent over $L_0$
will also be so near $L_0$. Dualising $\t E $, if
$\lambda_1,\ldots,\lambda_k \in \Gamma(L_0,\t E^*(a_r))$ are
independent, then they can be extended to sections of $\t
E^*(a_r)$  over a neighbourhood of $L_0$ in $\t U$ which are
linearly  independent at each point in this neighbourhood.
Rewriting $\t E\!\mid  _{L_0}$ in the form $\t E\!\mid_{L_0} \,=
\oplus^m_{i=1} V_i(b_i)$ where $V_i $ is a vector space with
$b_1<b_2<\ldots<b_m$, the last statement implies that the
projection $\t E\!\mid _{L_0} \to V(b_m)$ extends to an
epimorphism $\t E \to V(b_m)$ in a neighbourhood of $L_0$.
Using induction on rank, this gives the following local
description  of bundles on the blowup of a two-dimensional ball
at a point:

\bigskip
\noindent {\bf Proposition~2.5.}\quad {\sl Suppose $ \t
E\!\mid_{L_0} = \oplus^m_{i=1} V_i(b_i)$ where $V_i $ is a
$d_i$-dimensional vector space and $b_1<b_2<\cdots<b_m$.  Then in
a  neighbourhood of $L_0$ in $\t U$ there is a filtration
$$
0 =
F_0 \;\subset\; F_1\;\subset\; \ldots \;\subset F_m =\, \t E
\eqno (2.6)
$$
of $\t E$ by vector bundles $F_k $ such that
$F_k/F_{k-1} \simeq V_k(b_k) $ and such that $ F_k\!\mid_{L_0}\,
= \bigoplus^k_{i=1} V_i(b_i) $.  \qed }

\bigskip
\noindent {\bf Corollary~2.7.} \quad{\sl If  $ \t
E\!\mid _{L_0} = \oplus^r_{i=1} {\cal O}(a_i)$ with $|a_i-a_j|
\le 2$ for all $i,j$ then $\t E \simeq  \oplus^r_{i=1} {\cal
O}(a_i)$ in a neighbourhood of $L_0$.}

\medskip
\noindent{\bf Proof: } A bundle $F$ on $\t U$ which is
given as an extension $0 \to V(a) \to F \to W(b) \to 0$ for some
vector spaces $V,W$ is determined by an element of $H^1(\t
U,W^*\otimes V(a-b))$. If the extension splits on $\pi^{-1}(x_0)$
then this class lies in the image of $H^1(\t U,W^*\otimes
V(a-b+1)) \to H^1(\t U,W^*\otimes V(a-b))$, and  if $a-b \ge -2$
then the former group vanishes, by (2.1). \qed

\bigskip
The preceding discussion provides some insight into
well-known constructions of Serre [Ser] and of Schwarzenberger
[Sch], both described in [OSS].  A bundle $\t E$ on $\t X$ is the
pull-back of a bundle $E$ from $X$ if and only if $\t E$
restricts trivially to every component of the exceptional divisor
$\t Y$; if this is the case, then  necessarily $E = \pi_*\t E$.
Let $L_1,L_2$ be line bundles on $X$  and let $L_i$ also denote
$\pi^*L_i$ on $\t X$. Extensions of the  form $\,\t X: \, 0 \to
L_1(-1) \to \t E \to L_2(1) \to 0\,$ are  classified by $H^1(\t
X,L_1L_2^*(-2))$ which can be computed from  the Leray spectral
sequence for $\pi$ using (2.1). This gives  an exact sequence
$$
0 \to H^1(X,L_1L_2^*) \to H^1(\t X,L_1L_2^*(-2)) \to
H^0(Y,L_1L_2^*\otimes {\rm det}\;{\cal N}_Y) \to
H^2(X,L_1L_2^*)\eqno(2.8)
$$
and the bundle $\t E$ corresponding
to an element of $H^1(\t X, L_1L_2^*(-2))$ is trivial on $\pi^{-
1}(x_0)$ if and only if the element of $L_1L_2^*\otimes{\rm
det}\;{\cal N}_{Y,x_0}$ obtained from (2.8) is non-zero.   Using
(2.1) again, the bundle $E$ is given by an exact sequence $0 \to
L_1 \to E \to L_2\otimes {\cal I}_Y \to 0$, so $E\otimes L_1^*$
has a section vanishing precisely at $Y$ and ${\cal N}_{Y}$ has
been extended to the bundle $E^*\otimes L_2$.  If $X$ is compact,
the Chern classes of $E$ are given by $c_1(E) = c_1(L_1) -
c_1(L_2)$ and $c_2(E) = c_1(L_1)c_1(L_2) + PD([Y])$ where $PD[Y]$
denotes the Poincar\'e dual of $[Y]$.

More generally, this construction applies when $Y$ is an
arbitrary codimension $2$ locally complete intersection in a
complex manifold $X$ provided that ${\rm det}\;{\cal N}_Y$ can be
extended to a line bundle on $X$.  Indeed, with the exception of
Corollary~2.7, all of the results of this section so far
presented  remain valid if $X$ has arbitrary dimension and $Y$ is
a codimension two locally complete intersection; modifications to
the proofs above are straightforward.

\bigskip

If $X$ is compact the charge on $\t E$ can be estimated in terms
of its splitting on $L_0=\pi^{-1}(x_0)$ and the charge on its
direct image:

\bigskip

\noindent{\bf Proposition~2.9.} \quad {\sl Let $\blowup X$ be the
blowup of the compact surface $X$ at $x_0$, with $L_0 =
\pi^{-1}(x_0)$.  If $\t E$ is an $r$-bundle on $\t X$ such that
$\t E\!\!\mid_{L_0} = \oplus_{i=1}^r{\cal O}(a_i)$, then for $a:=
\sum_i a_i$ and $E := (\pi_*\t E)^{**}$ it follows
$$
C(E) +
{1\over 2}\sum_{i=1}^r|a_i-{a\over r}| + {\bar a\over 2r}(2r-\bar
a - n) \;\;\le\;\; C(\t E) \;\;\le \;\; C(E)+{1\over
2}\sum_{i=1}^r(a_i-{a\over r})^2  \eqno (2.10)
$$
where $\bar a
\equiv a \;({\rm mod}\; r)$, $0\le a <r$ and $n := \#\{a_i\mid
a_i\le a/r\}$.  Moreover, equality holds in the first case iff
$\t E\simeq (\pi^*E)(k)$ for some $k \in \m Z$ and in the second
iff $\t E \simeq \oplus {\cal O}(a_i)$ in a neighbourhood of
$L_0$.}

\medskip
\noindent{\bf Proof:} \quad  Since $\chi(\t
E)=\chi(\pi_*\t E)-\chi(\pi^1_*\t E)$ and $c_1(\t
X)=\pi^*c_1(X)+c_1({\cal O}(1))$ the Riemann-Roch formula gives
$$
C(\t E)=C(\pi_*\t E)+\chi(\pi^1_*\t E)-{a(a+r)\over 2r}=C(E)+
dim(E/\pi_*\t E)+\chi(\pi^1_*\t E)-{ a(a+r)\over 2r}\;.
\eqno(2.11)
$$
If $a_i\le 0$  for all $i$ and $\t E$ splits in a
neighbourhood of $L_0$ as a direct sum of line bundles and then
$\pi_*\t E$ is locally free and $\chi(\pi^1_*\t
E)=-\sum\chi({\cal  O}(a_i))+\sum\chi(\pi_*{\cal O}(a_i))
=(1/2)\sum a_i(a_i+1)$, giving $C(\t
E)=C(E)+(1/2)\sum(a_i-a/r)^2$ in this case.  If $\t E$ does not
split in a neighbourhood of $L_0$ then $\chi(\pi^1_*E)$ is
strictly less than  the corresponding number in the split
case---this follows easily by induction on rank using
Proposition~2.5.  Since (2.10) is invariant under tensoring $\t
E$ by ${\cal O}(-k)$, this proves the upper bound.

\smallskip

To obtain the lower bound, first twist $\t E$ by a suitable line
bundle so that $0 \le a <r$. Using Lemma~2.3 as in the proof of
Proposition~2.5, in a neighbourhood of $L_0$ the bundle $\t E$
can be written as an extension $0 \to A \to \t E \to  B \to 0$
which splits on $L_0$ and where $A\!\!\mid_{L_0} = \oplus_{a_i\le
0} {\cal O}(a_i)$ and $B\!\!\mid_{L_0}= \oplus_{a_i > 0} {\cal
O}(a_i)$.  Since the extension splits on $L_0$, it lies in the
image of $H^1(\t U,B^*\otimes A(1))$, so there is a bundle $\t
E_0$ on $\t U$ which is given by a compatible extension $0 \to
A(1) \to \t E_0 \to B \to 0$.  Compatibility of the extension
implies that there is an exact sequence $0 \to \t E_0 \to \t E
\to A\!\!\mid_{L_0} \to 0 $   giving $\pi_*\t E_0(-1) =\pi_*\t
E(-1)$ and $0 \to \pi^1_*\t E_0(-1) \to \pi^1_*\t E(-1) \to
\pi_*^1 A(-1)\!\!\mid_{L_0} \to 0$.   Hence $\chi(\pi^1_*\t
E(-1))\ge \chi(\pi_*^1 A(-1)\!\!\mid_{L_0}) = -\sum_{a_i \le
0}a_i$, so from (2.11) applied to $\t E(-1)$  it follows $C(\t
E)=C(\t E(-1))\ge C(E) + \sum_{a_i \le 0}|a_i| + a(r-a)/2r$. If
$a=0$ replace $\t E$ by $\t E^*$, otherwise replace $\t E$ by
$(\t E^*)(1)$ to obtain another lower bound $C(\t E) \ge
C(E)+\sum_{a_i \ge 0}|a_i| + a(r-a)/2r$. Averaging the two gives
$$
C(\t E) \;\ge \;{1\over 2}\sum_{i=1}^r|a_i| + {a(r-a)\over 2r}
\qquad\hbox{ if \quad $0\le a <r$.} \eqno(2.12)
$$
The expression
on the left of (2.10) follows immediately from this in the
general case by simple calculation.  Clearly equality holds in
(2.12) iff $a_i=0$ for all $i$. \quad \qed

\bigskip
\noindent The lower bound in (2.10) is sharp: given any
splitting type $\oplus{\cal O}(a_i)$ for a bundle $\t E$ on
$L_0$, an argument by induction on rank shows that there is a
bundle $\t E_0$ which splits minimally on (and hence in a
neighbourhood of) $L_0$ which is given as an extension $0 \to
A(1) \to \t E_0 \to B \to 0$ where $A$ and $B$ are bundles of the
form described above.  The inequalities (2.10) appear in  [FM2]
(Remark 5.4) in the case of rank $2$ bundle with $c_1=0$.

\bigskip
As before, let $X$ be a complex surface and let $\blowup
X$ be the blowup of $X$ at the point $x_0 \in X$.  The local
description provided by Proposition~2.5 can be combined with  a
holomorphic version of Taubes' ``cut-and-paste" construction [T]
to provide a global description of bundles on $\t X$.

Let $\t E$ be a holomorphic $r$-bundle on $\t X$, and let $U
\subset X$ be a neighbourhood of $x_0$ isomorphic to a ball. Then
$\t E$ can be viewed as comprised of two pieces,  namely the
bundle $(\pi_*\t E)^{**}$ on $X$ and the bundle $\t E\!\mid_{\t
U}$ on $\t U = \pi^{- 1}(U)$; the two pieces are glued  together
by means of the isomorphism $\pi_*\t E \simeq (\pi_*\t E)^{**}$
over $U\backslash \{x_0\} \simeq\t U \backslash L$.

Conversely, given $r$-bundles $E_0$ on $X$ and $E_1$ on $\t U$,
the two can be glued together by means of an isomorphism
$\rho\,\colon \pi_*E_1 \to E_0$ over $U \backslash\{x_0\}$,
extending uniquely to $U$ as an isomorphism $(\pi_*E_1)^{**} \to
E_0\!\!\mid\!_U$  to define a bundle $E_0\#_{\rho}E_1$ on $\t X$.
If $(E'_0,E'_1,\rho ')$ is another triple of such objects such
that ${E_0 \#_{\rho}E_1} \isoto {E'_0 \#_{\rho'}E'_1}$, then by
Hartogs' theorem the isomorphism $E_0\!\!\mid_{X \backslash \b U}
\;\isoto \; E'_0\!\!\mid_{X \backslash \b U}$ extends to an
isomorphism $\phi_0$ over $X$, and if $\phi_1$ denotes the
induced isomorphism $E_1 \to E'_1, $ it follows that $\rho' =
\phi_0^{}\, \rho \,\pi_*\phi_1^{-1}$.

If $\rho_0 \colon (\pi_*E_1)^{**} \to E_0\!\!\mid_U$ is fixed and
$\rho$ is any other such isomorphism, then $\rho\rho_0^{-1}$ is
an automorphism of $E_0$ over $U$.  Hence the following
description is obtained:

\bigskip
\noindent{\bf Proposition~2.13.} \quad {\sl Let $E_0$ be
a bundle on $X$, $E_1$ be a bundle on a neighbourhood
$\pi^{-1}(U)$ of $\pi^{-1}(x_0)$ in the blowup $\t X$ of $X$ at
$x_0$ and $\rho_0 \colon (\pi_*E_1)^{**} \isoto E_0\!\!\mid_U$ be
given. Then isomorphism classes of vector bundles $\t E$ on $\t
X$ such that $(\pi_*\t E)^{**} \simeq E_0$ and $\t E \simeq E_1$
in a neighbourhood of $\pi^{-1}(x_0)$ are parameterised by the
stalk of the skyscraper sheaf ${\cal
A}ut(E_0)/\rho_0^{}\pi_*{\cal A}ut(E_1)\rho_0^{-1}$ at $x_0$,
modulo the left action of $\Gamma(X,Aut(E_0))$. \quad \qed}

\bigskip
\noindent{\bf Remarks:} \quad {\bf 1.}\quad If $E_1 =
\bigoplus\,{\cal O}(a_i)$, the space ${\cal
A}ut((\pi_*E_1)^{**})/ \pi_*{\cal A}ut(E_1)$ is easily identified
with the total space of some homogeneous vector bundle over a
flag  manifold.  Two simple examples which will be of some
relevance subsequently are the cases $\,E_1 = {\cal O}(1) \oplus
{\cal O}^{r- 1}\,$ and $\,E_1 = {\cal O}(-1)\oplus {\cal O}(1)
\oplus {\cal O}^{r-2}$, for which the corresponding skyscrapers
are respectively ${\m{F}}_1({\m{C}}^r)$ and the total space of
the bundle $2{\cal O}(1,1)$ over ${\m{F}}_{1,r-1}({\m{C}}^r)$;
(the zero section corresponds to those bundles which extend to
$\t{\m P}_2$ as a direct sum of line bundles).

\medskip
\noindent{\bf 2.}\quad This description provides a
simple way to construct bundles on $\t X$ from bundles on $X$,
but in contrast with the construction of Serre/Schwarzenberger
the bundles produced this way are all non-trivial on the
exceptional divisor.  However, if $\sum_i a_i = 0$ the generic
deformation of the bundle $\t E$ (or $E_1$) will be trivial on
$L$, and the earlier construction can be seen as a deformation of
a bundle restricting to ${\cal O}(-1) \oplus {\cal O}(1)$ on $L$;
(the existence of such deformations is discussed further in \S4).

\bigskip
An effective description of the spaces of holomorphic
bundles on $\t X$ requires such a description for the spaces of
bundles in a neighbourhood of $L_0$, but that given by
Proposition~2.5 has some redundancy: the filtration (2.6) is not
uniquely determined. However, by gluing a bundle on a
neighbourhood of $L_0$ to the trivial bundle on $\m P_2$ the
classification problem becomes that of determining the bundles on
$\t{\m P}_2$ which are trivial in a neighbourhood of the line
$L_{\infty}$ at infinity. A trivialisation of such a bundle in a
neighbourhood of $L_{\infty}$ is determined by its restriction to
$L_{\infty}$, so isomorphism  classes of pairs $(E_1,\rho)$ where
$\rho$ is a trivialisation of $(\pi_*E_1)^{**}$ in a
neighbourhood of $x_0$ correspond to isomorphism classes of
pairs $(\t E,\varphi)$ where $\t E$ is a bundle on $\t{\m{P}}_2$
such  that $(\pi_*\t E)^{**}$ is trivial and $\varphi$ is a
trivialisation of $\t E$ on $L_{\infty}$.

The space $\t{\m P}_2$ is isomorphic to the Hirzebruch surface
$H_1$ and bundles on this space have been studied in [B1].  Using
the lemma of \S1 of that reference, a monad description of all
holomorphic bundles on $H_1$ trivial on $L_{\infty}$ is easily
given, and the precise condition on such monads for the
corresponding bundles to be trivial in a neighbourhood of
$L_{\infty}$ is easily calculated; ([B5]).

\bigskip
\bigskip
\bigskip

\noindent {\bf 3.\quad Structure of moduli spaces I.}

\bigskip
The ``cut-and-paste" construction of the previous
section makes no reference to questions of stability, an issue
which is considered in this section.  When the discussion is not
limited to a single bundle but rather to whole moduli spaces, the
results generally apply only in the case that $b_1(X)$ is even;
the reason for this can be traced to the failure of Lemma~1.8
when $b_1(X)$ is odd.

\medskip
The definition of stability requires a hermitian metric,
and throughout this section such a metric (positive (1,1)-form)
$\omega$ is a fixed on $X$.  The metrics to be used on blowups of
$X$ are the same as those used in [B2], the construction of which
will be briefly recalled here for convenience.

Let $\blowup X$ be the blowup of $X$ at $x_0 \in X$.  Let $L :=
\pi^{-1}(x_0)$ be the exceptional divisor so $\pi^*\omega$ is
everywhere non-negative and is degenerate only in directions
tangent to $L$.  Let $\sigma$ be $i/2\pi$ times the curvature
form of any hermitian connection on the line bundle ${\cal O}(-L)
=: {\cal O}(1)$ restricting positively to $L$, and let
$\omega_{\epsilon} := \pi^*\omega + \epsilon \sigma$ for
$\epsilon > 0$; (recall ${\cal O}(L)\!\mid_{L} \simeq {\cal
O}_L(-1)$).  It follows that if $\epsilon$ is sufficiently small
then $\omega_{\epsilon}$ defines a positive form in a
neighbourhood of $L$; if $\sigma$ is compactly supported in $\t
X$ then $\omega_{\epsilon}$ is everywhere positive for
sufficiently small $\epsilon$.  If $\omega$ is $\dbd$-closed and
$\sigma$ is compactly supported, it follows from the fact that
$L$ has self-intersection $-1$ that $Vol(\t X,\omega_{\epsilon})
= Vol(X,\omega) -\epsilon^2/2=V-\epsilon^2/2$, and if $\omega$ is
$d$-closed, then so too is $\omega_{\epsilon}$.

A useful model to keep in mind is the following: if $x_0$
corresponds to the origin in local  holomorphic coordinates
$\{z^a\}$,  the orientation-reversing map  $z^a \mapsto
z^a/|z|^2$ lifts to the blowup to define an  isomorphism of a
neighbourhood of $L$ with a neighbourhood of a line  in
$\m{CP}_2$, realising $\t X$ as the connected sum $\t X
\simeq_{diffeo} X\#\b{\m{CP}}_2$.   Under this diffeomorphism,
the pullback of $\omega_1 = (i/2)\ddb(|z|^2 + \log|z|^2)$  is
conformal to the Fubini-Study metric.  The form $\sigma$ can be
taken to be $(i/2)\ddb\log(\psi(|z|^2))$ where $\psi(t)$ is a
smooth function which is the identity near $0$ and a positive
constant for $t \ge t_0$.  Pulling back under the ``dilations" $z
\mapsto \epsilon^{-1/2} z$ and rescaling by $\epsilon$ gives the
metric $\omega_{\epsilon}$, ``stretching out" the neck of the
connected sum as in [D3].

Now let $\blowup X$ be a modification of $X$ consisting of $n$
successive blowups (at simple points), and let $\sigma_i$ be a
closed smooth (1,1)-form on $\t X$ corresponding as in the last
paragraph to the $i$-th blowup.  Let $\m R^n_+ := \{ \alpha =
(\alpha_1,\ldots,\alpha_n) \in {\m{R}}^n\!\mid \alpha_i >0 \,,\;
i = 1,\ldots,n \}$ and for $\alpha \in \m R^n_+\,$ let $
\rho_{\alpha} := \sum \alpha_i\sigma_i$, so
$\rho_{\alpha}\!\cdot\!\rho_{\alpha} = - \sum\alpha_i^2 =:
-|\alpha|^2$ and $ \omega_{\alpha} := \pi^*\omega +
\rho_{\alpha}$ is positive for $|\alpha|$ sufficiently small;
(this definition differs slightly from that in [B2] where $\rho$
has the opposite sign).  A vector $\alpha \in \m R^n_+$ is called
{\it suitable\/} if $\omega_{\alpha}$ is a positive form on $\t
X$.

\bigskip
Let $\blowup X$ be a blowup of the compact surface $X$,
equipped with a metric of the form $\omega_{\alpha}$ as above.
Let $\t E$ be a holomorphic bundle on $\t X$, and let $\t A
\subset \t E$ be a bundle of rank $a$ included in $\t E$ as a
subsheaf.  By definition of $\nu_{\bullet}(*)$,
$$
\nu_{\t E}(\t
A,\omega_{\alpha})\; = \;\nu_{\pi_*\t E}(\pi_*\t A) +
\rho_{\alpha}\!\cdot\![{a\over r}c_1(\t E) - c_1(\t
A)]\;.\eqno(3.1)
$$
If $\t A$ has torsion-free quotient $\t C$ of
rank $c>0$ then (1.2) and (1.4) give
$$
\eqalignno{C(\t E) \; =\;
& C(\t A) + C(\t C)&(3.2)\cr
&+ {r\over 2ac}\Big[\,\norm{{a\over
r}c_1(\pi_*\t E) - c_1(\pi_*\t A)}^2_{\omega} + \norm{ {a\over
r}c_1(\t E)-c_1(\t A)}_Q^2 - V^{- 1}\nu_ {\pi_*\t E} (\pi_*\t
A)^2\,\Big]\,,}
$$
where $\norm{x}^2_Q := - (x
-\pi^*\pi_*x)\!\cdot\! (x - \pi^*\pi_*x)$ for $x \in H^2(\t
X,\m{Q})$.

Suppose now that $\t E$ is semi-stable but not stable with
respect to $\omega_{\alpha}$, and that $\t A$ destabilises $\t
E$.  Then from (3.1) it follows that $\nu_{\pi_*\t E}(\pi_*\t A)
= -\rho_{\alpha}\!\cdot\![{a\over r}c_1(\t E) - c_1(\t A)]$ so
(3.2) implies
$$
\eqalignno{C(\t E) \;\ge \; & C(\t A) + C(\t
C)&(3.3)\cr
& + {r\over 2ac}\Big[\,\norm{{a\over r}c_1(\pi_*\t E)
- c_1(\pi_*\t A)}^2_{\omega} + (1-V^{-1} |\alpha|^2)\norm{
{a\over r}c_1(\t E) - c_1(\t A)}_Q^2\,\Big]\,.}
$$
Since $\t A$
destabilises $\t E$ for the metric $\omega_{\alpha}$, it follows
from the semi-stability of $\t E$ that both $\t A$ and $\t C$ are
also semi-stable (with respect to this metric), implying $C(\t
A), C(\t C)$ are non-negative.  Hence (3.3) yields a uniform
bound on $\norm{(a/r)c_1(\t E) - c_1(\t A)}$ which involves only
$C(\t E)$ and $r$ (if $|\alpha|$ is suitably bounded from above).

\medskip
With these preparations in hand, the following result
summarises most of the important relationships between stability
on $X$ and $\t X$.  Parts of it appear in [FM2] (Theorem~5.5) and
[Br] (Theorem~4) in the case of bundles of rank $2$.

\bigskip
\noindent {\bf Proposition~3.4.} \quad {\sl Let $\t E$
be an $r$-bundle on $\t X$.
\smallskip
\item{(a)~} If $\t E$ is
$\omega_{\epsilon\alpha}$-stable for all sufficiently small
$\epsilon>0$ then $\pi_*\t E$ is $\omega$-semi-stable;
\item{(b)~} If $\t E=\pi^*E$ for some bundle $E$ on $X$ and $\t
E$ is $\omega_{\alpha}$-stable, then $E$ is $\omega$-stable;
\item{(c)~} If $\t E$ is $\omega_{\alpha}$-(semi-)stable and
$\pi_*\t E$ is $\omega$-semi-stable, then $\t E$ is
$\omega_{\epsilon\alpha}$-(semi-)stable for all $\epsilon \in
(0,1]$;
\item{(d)~} If   $\pi_*\t E$ is stable, then $\t E$ is
$\omega_{\alpha}$-stable for all suitable $\alpha \in \m R^n_+$
sufficiently small.  In fact, if $\nu_{\pi_*\t E}(A) \ge \delta
>0$ for all  $A \subset \pi_*\t E$ with non-zero torsion-free
quotient, then $\t E$ is $\omega_{\alpha}$-stable for all
suitable $\alpha \in \m R^n_+$ such that $|\alpha| <
\delta\sqrt{{2V\over 2\delta^2+rVC(\t E)}}$.  \par}

\medskip

\noindent{\bf Proof: }   (a) \quad By Lemma~2.3, after twisting
$\t E$ with a suitable line bundle it can be supposed that
$\pi_*\t E$ is locally free. If $A \subset \pi_*\t E$ has
torsion-free quotient, then $\pi^*A \subset \pi^*\pi_*\t E
\subset \t E$.  If $\t A$ is the maximal normal extension of
$\pi^*A$ in $\t E$ then $\mu(\pi^*A,\omega_{\alpha}) \le \mu(\t
A,\omega_{\alpha})$.  Replacing $\alpha$ by $\epsilon\alpha$ in
(3.1) and letting $\epsilon \to 0$ gives $\nu_{\pi_*\t E}(A)\ge
0$.

\noindent (b) \quad If $A\subset E$ has torsion-free quotient and
$\t A\subset \t E$ is the maximal normal extension of $\pi^*A$ in
$\t E=\pi^*E$, then $\mu(A,\omega)=\mu(\pi^*A,\omega_{\alpha})\le
\mu(\t A,\omega_{\alpha})<\mu(\t
E,\omega_{\alpha}=\mu(E,\omega)$.

\noindent (c) \quad Suppose $\t A \subset \t E$ has torsion-free
quotient. If $\nu_{\t E}(\t A,\omega_{\epsilon\alpha}) <0$ for
some $\epsilon \in (0,1)$ then since $\nu_{\pi_*\t E}(\pi_*\t A)
\ge 0 $ by hypothesis it would follow that
$\epsilon\rho_{\alpha}\!\cdot\![{a\over r}c_1(\t E) - c_1(\t A)]
<0$; this would imply $\nu_{\t E}(\t A,\omega_{\alpha}) <0$ also,
contradicting the hypotheses.  If $\t A$ destabilises $\t E$ for
$\omega_{\epsilon\alpha}$ then both $\nu_{\pi_*\t E}(\pi_*\t A)$
and $\rho_{\alpha}\!\cdot\![{a\over r}c_1(\t E) - c_1(\t A)]$
must be zero. Otherwise, one must be strictly positive and
therefore $\t E$ must be strictly stable with respect to
$\omega_{\epsilon\alpha}$ for all $\epsilon \in (0,1]$.

\noindent  (d) \quad If $\t A \subset \t E$ has torsion-free
quotient, the first Chern class of $\t A$ restricted to any
irreducible component of the exceptional divisor is bounded above
by a constant depending on the maximum of the first Chern classes
of the line bundles in the decomposition of $\t E$ on that
component.  The proof of Lemma~5 of [B2] now applies to show that
if $\t E$ has at least one non-trivial subsheaf, then there
exists such a subsheaf $\t A$ with torsion-free quotient $\t C$
which maximizes $\mu(-,\omega_{\epsilon\alpha})$ for any
sufficiently small $\epsilon > 0$.  It then follows from (3.1)
that $\t E$ is stable with respect to $\omega_{\epsilon\alpha}$
for all $\epsilon$ sufficiently small.

If $\alpha$ is suitable and satisfies the inequality of (d) and
$\t E$ is not $\omega_{\alpha}$-stable, then since stability is
an open condition on the metric there exists $\epsilon \in (0,1]$
such that $\t E$ is $\omega_{\epsilon\alpha}$-semi-stable but not
stable.  Then $\delta \le \nu_{\pi_*\t E}(\t A) =
-\epsilon\rho_{\alpha}\cdot((a/r)c_1(\t E)-c_1(\t A)) \le
|\alpha| \norm{c_1(\t E)-c_1(\t A)}_Q$.  Replacing $\alpha$ by
$\epsilon\alpha$ in (3.3) gives the bound $\norm{c_1(\t E)-c_1(\t
A)}_Q^2 \le rVC(\t E)/(2V-2\epsilon^2|\alpha|^2)$, which gives
the contradiction $\delta < \delta$ after a little algebra. \quad
\qed

\bigskip

A simple corollary of the proposition is that if $\t E =
E_0\#_{\rho}E_1$ is a bundle on $\t X$ constructed by the gluing
construction of the last section and if $E_0$ on $X$ is
$\omega$-stable, then $\t E$ is $\omega_{\alpha}$-stable for all
suitable $\alpha \in \m R^n_+$ sufficiently small.  The gluing
construction, combined with the existence of Hermitian-Einstein
connections on stable bundles, can thus be viewed as a
holomorphic interpretation of Donaldson's ``connected sums of
connections" theorem [D3] where one of the summands is (a
connected sum of copies of) $\b{\m{P}}_2$.

\medskip
By Lemma~1.8,  Proposition~2.9 and part (d) of
Proposition~3.4, when $b_1(X)$ is even it follows that for any
bundle $\t E$ on $\t X$ of rank $r$ and charge $\le C_0$ such
that $\pi_*\t E$ is stable,  $\t E$ is $\omega_{\alpha}$-stable
for all suitable $\alpha \in \m R^n_+$ satisfying a uniform bound
independent of $\t E$.  In particular, the pull-backs  from $X$
of stable bundles of bounded ranks and charge are all stable with
respect to the same metrics on the blowup.  The restriction that
$b_1(X)$ be even gives the following strengthening of
Proposition~3.4:

\bigskip

\noindent{\bf Proposition~3.5.} \quad { \sl Suppose that $b_1(X)$
is even. For any $r_0, C_0 >0$ there exists $\epsilon_0 =
\epsilon_0(r_0,C_0,\omega)$ with the property that any bundle on
a blowup $\t X$ of $X$ of rank $\le r_0$ and charge $\le C_0$
which is stable with respect to $\omega_{\alpha}$ for some
suitable $\alpha \in \m R^n_+$ satisfying $|\alpha_0|
<\epsilon_0$ is stable with respect to $\omega_{\epsilon\alpha}$
for all $\epsilon \in (0,\epsilon_0/|\alpha|)$.   Moreover, any
bundle which is semi-stable with respect to $\omega_{\alpha}$ and
has semi-stable direct image is semi-stable with respect to
$\omega_{\epsilon\alpha}$ for all $\epsilon \in
(0,\epsilon_0/|\alpha|)$.}

\medskip
\noindent {\bf Proof:} \quad Set $\epsilon_1 :=
\delta_0\sqrt{2V/(2\delta_0^2+rVC(\t E))}$ where $\delta_0$ is
as in Lemma~1.8.  Suppose that there exists suitable $\alpha_1
\in \m R^n_+$ with $|\alpha_1| < \epsilon_1$ and a bundle $\t
E_1$ of rank $r_1\le r_0$ on $\t X$ with $C(\t E_1) \le C_0$
which is $\omega_{\alpha_1}$-stable  but which is not
$\omega_{\delta_1\alpha_1}$-stable for some $\delta_1 \in
(0,\epsilon_1/|\alpha_1|)$, where $\delta_1\alpha_1$ is suitable.
Since stability is an open condition on the metric, by altering
$\delta_1$ if necessary it can be supposed that $\t E_1$ is
semi-stable but not stable with respect to
$\omega_{\delta_1\alpha}$, and hence that there exists $\t A_1
\subset \t E_1$ destabilising as in the discussion preceding
Proposition~3.4.  Set $\epsilon_2 :=
(1/2)\,min\{|\alpha_1|,\delta_1|\alpha_1|\}$  and repeat this
procedure, generating a sequence $\{(\epsilon_i, \alpha_i,\t
E_i,\t A_i,\delta_1)\}$ by iteration.  By construction,
$$
\eqalignno{\nu_{\pi_*\t E_i}(\pi_*\t
A_i)+\rho_{\alpha_i}\!\cdot\! [(a_i/r_i)c_1(\t E_i) - c_1(\t
A_i)]\;&>\;0&(3.6)(a)\cr
\nu_{\pi_* \t E_i}(\pi_*\t A_i)
+\delta_i\rho_{\alpha_i} \!\cdot\! [(a_i/r_i)c_1(\t E_i) - c_1(\t
A_i)] \;&=\; 0\;,&(b)\cr}
$$
with $\delta_1|\alpha_{i}| <
\epsilon_{i}<\epsilon_{i-1}/2$.  Since $1 \le a_i \le r_i \le
r_0$ and $C(\t E_i) \le C_0$ the uniform bound on the norms
$\norm{(a_i/r_i)c_1(\t E_i) - c_1(\t A_i)}$ provided by (3.3)
implies that there is a subsequence for which $(a_i/r_i)c_1(\t
E_i)-c_1(\t A_i)$ is constant.  Since $b_1(X)$ is even,
$\nu_{\bullet}(*)$ is topological and therefore constant on this
subsequence.  If the subsequence is infinite, then $\epsilon_i
\to 0$, so (3.6)(b) implies $\nu_{\pi_* \t E_i}(\pi_*\t A_i)$ is
eventually $0$ on this subsequence and therefore so too is
$\rho_{\alpha_i} \!\cdot\![(a_i/r_i)c_1(\t E_i) - c_1(\t A_i)]$;
this however contradicts (3.6)(a).  Thus the subsequence must be
finite, implying the original sequence terminated, which in turn
implies the first statement of the proposition.

\smallskip
To prove the second statement,  suppose  that $\t E$
is $\omega_{\alpha}$-semi-stable  for some suitable $\alpha \in
\m R^n_+$ with $|\alpha|<\epsilon_0$, and that $\pi_*\t E$ is
semi-stable.  If $\t E$ is $\omega_{\alpha}$-stable then the
first part of the proposition applies, so it can be assumed that
$\t E$ is not stable.  If $\t A \subset \t E$ has torsion-free
quotient and destabilises $\t E$ with respect to
$\omega_{\alpha}$, then $\nu_{\pi_*\t E}(\pi_*\t A)=
-\rho_{\alpha}\cdot [(a/r)c_1(\t E)-c_1(\t A)] \le |\alpha|\,
\norm{(a/r)c_1(\t E)-c_1(\t A)}_Q < \delta_0$ by construction of
$\epsilon_0$ and by (3.3), so by Lemma~1.8 it follows
$\nu_{\pi_*\t E}(\pi_*\t A)=0=\rho_{\alpha}\cdot [(a/r)c_1(\t
E)-c_1(\t A)]$.  Hence $\nu_{\t E}(\t
A,\omega_{\epsilon\alpha})=0$ for all $\epsilon$. \quad \qed

\bigskip
\noindent In general, it is not the case that the moduli
spaces of $\omega_{\alpha}$-stable holomorphic structures  on a
given topological bundle over $\t X$ are  independent of  $\alpha
\in \m R^n_+$ once $|\alpha|$ is sufficiently small. If $\t E$ is
a bundle on $\t X$  of rank $ \le r_0 $  and charge $ \le C_0$
which is stable with  respect to $\omega_{\alpha}$ but not stable
with respect to $\omega_{\beta} $ for $|\alpha| ,\, |\beta| \le
\epsilon_0$, then it follows easily as in the proof of
Proposition~3.5 that there exists $c \in H^2(X,\m Z)^{\perp}
\subset H^2(\t X,{\m{Z}}) \cap H^{1,1}(\t X)$ with $\norm{c}_Q
\le \sqrt{r_0^3VC_0/(2V-2\epsilon_0^2)}$ such that
$\rho_{\alpha}\cdot c > 0$ and $\rho_{\beta}\cdot  c \le 0$,
namely $c= (ac_1(\t E)-rc_1(\t A))-(a\pi^*c_1(\pi_*\t
E)-r\pi^*c_1(\pi_*\t A))$  for some destabilising $\t A\subset \t
E$.  The moduli spaces will be independent of suitable $\alpha
\in \m R^n_+$ satisfying $|\alpha| < \epsilon_0$ provided that
$\alpha$ remains within one of the finitely many chambers of $\m
R^n_+$ cut out by the equations $\rho_{\alpha}\!\cdot c = 0 $ for
$c \in H^2(X,\m Z)^{\perp}$ with $\norm{c}_Q\le
\sqrt{r_0^3VC_0/(2V-2\epsilon_0^2)}$. Such a ``chamber structure"
for moduli spaces is quite well-known---see, e.g., [D4], [K2].

Moduli spaces also depend non-trivially on $|\alpha|$ in general:
it is not hard to construct examples of bundles on a blowup $\t
X$ which are stable with respect to $\omega_{\alpha}$ for some
suitable $\alpha$,  but which fail to be stable with respect to
$\omega_{\epsilon\alpha}$ for some $\epsilon\in (0,1)$.

\bigskip
\bigskip
\bigskip

\noindent {\bf 4.\quad Stabilisation and desingularisation.}

\bigskip
The appearance of sheaves and bundles which are
semi-stable but not stable represents a divergence between the
real analytical and the complex analytical descriptions of
moduli: whereas isomorphism classes of stable bundles and
irreducible Hermitian-Einstein connections are in one-to-one
correspondence, this fails to be true as soon as {\it stable} is
replaced by {\it semi-stable}: for example, if $A$ and $C$ are
stable bundles and $\mu(A) = \mu(C)$, then any extension of the
form $0 \to A \to E \to B \to 0$ defines a bundle $E$ with
$\mu(E) = \mu(C)$ and which is always semi-stable but not stable.
The bundle $E$ admits a Hermitian-Einstein connection if and only
if the extension splits.
In this section, the gluing construction of \S2 is used to
provide a mechanism for ``stabilising" a semi-stable bundle (or
torsion-free sheaf). The methods, which are strictly
sheaf-theoretical,  can also be used to by-pass some of the
technical difficulties encountered in [D5] to show that moduli
spaces of stable bundles of sufficiently large charge on a blowup
of a surface have open subsets of which are smooth;  this is
indicated in the second half of the section.  Finally, the same
methods are used to show that bundles on a blowup which are
topologically trivial on the exceptional divisor can be
approximated (off a finite set) by bundles which are
holomorphically trivial on the divisor.

\medskip

In general, a semi-stable sheaf ${\cal S}$ on the compact surface
$X$ determines a semi-stable bundle $\Sigma ({\cal S})$ admitting
a Hermitian-Einstein connection as follows: $\Sigma ({\cal S}) :=
\Sigma({\cal S}^{**})$ and if ${\cal A} \subset {\cal S} $ has
$\mu({\cal A}) = \mu({\cal S})$, then $\Sigma({\cal S}) :=
\Sigma({\cal A}) \oplus \Sigma({\cal S}/{\cal A})$.  It is
straightforward to verify by induction on rank that this
prescription is well-defined and uniquely determines the bundle
$\Sigma({\cal S})$.  This bundle has the same rank and
determinant as ${\cal S}$ and never has greater charge; it is a
direct sum of stable bundles all of the same slope (i.e., is
quasi-stable), and there are non-zero holomorphic maps ${\cal
S}\to \Sigma({\cal S})$, $\Sigma({\cal S}) \to {\cal S}^{**}$.

It is also convenient to introduce the notation $B(E)$ for a
semi-stable bundle $E$ to denote the set of points $x\in X$ for
which there is a semi-stable bundle $A$ with $\mu(A)=\mu(E)$ and
a sheaf inclusion $A\to E$ such that $A_x \to E_x$ {\it not\/} of
maximal rank; (note that the quotient $E/A$ must be torsion-free
else semi-stability of $E$ will be violated by the maximal normal
extension of $A$ in $E$). Again, it is easily verified by
induction on the rank of $E$ that $B(E)$ is finite.

\medskip
Let $E$ be a semi-stable $r$-bundle on $X$, and let $A$
be a semi-stable $a$-bundle with $\mu (A) = \mu (E)$ for which
there is a map $A \to E$ inducing a sheaf inclusion.   Pick a
point $x_0 \in X \backslash B(E)$ and let $\blowup X$ be the
blowup of $X$ at $x_0$.   If $E_1$ is any $r$-bundle on a
neighbourhood of $L := \pi^{-1}(x_0)$ and $\rho \colon \,E \to
(\pi_*E_1)^{**}$ is an isomorphism over this neighbourhood, then
if $E$ is in fact stable it follows from Proposition~3.4 that the
bundle $\t E := E \#_{\rho}E_1$ is stable with respect to
$\omega_{\epsilon}$ for all sufficiently small $\epsilon > 0$;
the more delicate and interesting case is when $E$ is not stable
which is henceforth assumed.

Now take $E_1 = {\cal O}(1) \oplus  {\cal O}^{r-1}$.  Up to
isomorphism, the bundle $\t E$ is determined by a non-zero
element, $\varphi$ say, of the vector space $E^*_{x_0}$, with any
non-zero multiple giving an isomorphic bundle: the correspondence
is given explicitly by taking direct images of the sequence $0
\to \t E \to \t E(-1) \to \t E(-1)\!\!\mid_L \to 0$ and using the
fact that $\pi_* \t E(-1)$ is locally free, hence equal to $E$.

Let $\t A \subset \t E$ have torsion-free quotient $\t C$ and
satisfy $\mu (\pi_*\t A) = \mu (E)$, and consider the commutative
diagram
$$
\matrix{0 &\rightarrow &\pi_*\t A &\rightarrow
&\pi_*\t E & \rightarrow & \pi_*\t C & \rightarrow & 0\quad
\hphantom{.} \cr
&&\downarrow && \downarrow {\vrule height 12 pt
depth 5 pt width 0 pt} && \downarrow \cr
0 &\rightarrow &
{}~(\pi_*\t A)^{**} & \rightarrow &~~(\pi_*\t E)^{**}
&\rightarrow &~~(\pi_*\t C)^{**} & \rightarrow & 0 \quad .}
$$
Since $\pi_*^1\t E = 0$, the cokernel of $\pi_*\t E \to \pi_*\t
C$ is $\pi_*^1\t A$.  Since $x_0 \notin B(E)$, $(\pi_*\t A)^{**}
\subset (\pi_*\t E)^{**} = E$ is a sub-bundle at $x_0$ and the
lower row must be exact near there, implying that $(\pi_*\t
A)^{**}/\pi_*\t A \to (\pi_*\t E)^{**}/\pi_*\t E $ is injective
and that $(\pi_*\t E)^{**}/\pi_*\t E \to (\pi_*\t C)^{**}/\pi_*\t
C$ is surjective.

If the kernel of $ \varphi \,\colon E_{x_0} \onto {\m{C}}$ does
not contain the image of $(\pi_*\t A)^{**}$, then the restriction
of $\varphi$ to $(\pi_*\t A)^{**}$ is non-zero and since
$(\pi_*\t E)^{**}/\pi_*\t E = {\m{C}}$, this implies that
$(\pi_*\t A)^{**}/\pi_*\t A = {\m{C}}$ also.  It follows in this
case that $\pi_*^1\t A = 0$ and $\pi_*\t C = (\pi_*\t C)^{**}$
implying (since $\pi_*^1\t C = 0$) that $\t C$ is locally free
near $L$, and by Lemma~2.3, that if $\t C\!\!\mid_L = \sum {\cal
O}(c_i)$ then $\sum_i c_i \le 0$.  Since $E_1\!\!\mid_L \to \t
C\!\!\mid_L $ is onto, all $c_i$ must be non-negative, so $c_i =
0$ for all $i$, and $\t A\!\!\mid_L = {\cal O}(1) \oplus \sum
{\cal O}$.

On the other hand, if the kernel of $\varphi$ {\it does\/}
contain the image of $(\pi_*\t A)^{**}$ then $\pi_*\t A =
(\pi_*\t A)^{**}$ and there is an exact sequence $0\to \pi^1_*\t
A \to \m C\to (\pi_* \t C)^{**}/\pi_*\t C \to 0$.  Since $\pi_*
\t A$ is locally free, Lemma~2.3 implies that near $L$, $\det \t
A= {\cal O}(a)$ for some $a\le 0$, and since $\det \t E= {\cal
O}(1)$, the same lemma implies that $\pi_*(\t C^{**})$ cannot be
locally free near $x_0$.  In the exact sequence $0\to \pi_*(\t
C^{**}/\t C) \to (\pi_*\t C)^{**}/\pi_*\t C \to (\pi_*(\t
C^{**}))^{**}/\pi_*(\t C^{**}) \to 0$ the last term is therefore
non-zero, implying the same of the middle term.  Thus from the
previous exact sequence it now follows that $\m C = \pi_*(\t
C)^{**}/\pi_*\t C = (\pi_*(\t C^{**}))^{**}/\pi_*(\t C^{**})$
implying $\t C$ is locally free, and that $\pi_*^1\t A=0$.

Since $\t A \to \t E$ is now a bundle map, it follows that if $\t
A$ splits on $L$ as $\sum {\cal O}(a_i)$, then all $a_i$ must be
$\le 1$.  Since $\pi_*^1\t A =0$, all $a_i$ must also be $\ge
-1$, so by Corollary~2.7 it follows that $\t A$ splits as $\sum
{\cal O}(a_i)$ in a neighbourhood of $L$.  Since $\pi_*\t A$ is
locally free, it must therefore be the case that all $a_i$ are in
fact $\le 0$.  Now replace the top row of the diagram above with
the exact sequence $0 \to \pi_*\t A(- 1)\to\pi_*\t
E(-1)\to\pi_*\t C(-1)$, the cokernel of the last arrow being
$\pi^1_*\t A(-1)$.  Since the middle term is now locally free, it
follows that the kernel of $\pi_*\t C(-1)\to \pi_*^1\t A(- 1)$ is
isomorphic to $(\pi_*\t C(-1))^{**}$, so $\pi_*\t C(-1)$ is
locally free and $\pi_*^1\t A(-1)=0$.  Consequently $\t A$ is
trivial and $\t C= {\cal O}(1)\oplus \sum {\cal O}$ near $L$ in
this case.

Thus the quotient $\t C$ is always locally free near $L$ and $\t
A\! \!\mid_L = \sum {\cal O}$ or $ {\cal O}(1) \oplus \sum {\cal
O}$ according to whether the kernel of $\varphi$ does or does not
contain the image of $(\pi_*\t A)^{**}$ respectively.  If it does
not, then $\mu (\t A,\omega_{\alpha}) = \mu (\pi_*\t A) -\alpha
/a = \mu (E) -\alpha /a <\mu (E) -\alpha/r = \mu (\t E
,\omega_{\alpha})$, so $\t E$ would be stable for all $\alpha >0$
sufficiently small if every such $\t A$ could be guaranteed to
fall into the second category.

In general, it is not be possible guarantee that $(\pi_*\t
A)^{**}$ should not be contained in the kernel of $\varphi$.
However, the following somewhat technical lemma shows that if the
construction is repeated at a number of points, a uniform upper
bound on the number of such ``bad" points can be given:
\bigskip
\noindent {\bf Lemma~4.1.} \quad {\sl Suppose $x_1,\dots,x_n \in
X\backslash B(E)$.  Then there are linear maps $\varphi_i
\,\colon E_{x_i} \to {\m{C}}, \; i = 1,\ldots,n$ with the
following property: for any stable bundle $A$ on $X$ of rank $a$
with $\mu (A) = \mu (E)$ and any non-zero map $A \to E$, the
number of maps $\varphi_i$ for which $A_{x_i} \subset
ker\,\varphi_i$ is at most $\,m-1$, where the multiplicity of $A$
in $\Sigma(E)$ is $\,ma -c\,, \; 0 \le c < a$. }

\medskip
\noindent {\bf Proof: } \quad By induction on the rank
$r$ of $E$.  Without loss of generality, there is a semi-stable
bundle $K$ with $\mu (K) = \mu (E)$ and a map $K \to E$ such that
the quotient $B$ is a torsion-free stable sheaf of rank $0 < b <
r$.  By the inductive hypotheses, maps $\varphi_i^K$ with the
requisite properties exist for $K$.

Suppose first that the extension $0 \to K \to E \to B \to 0 $ is
non-trivial.  If $A$ is a stable bundle of rank $a$ with $\mu (A)
=\mu (E)$ and there is a non-zero map $A \to E$, then by
stability of $A$ and $B$ the composition $A \to B$ is either $0$
or an isomorphism, but the latter is ruled out by the assumption
that the extension does not split.  Thus any such $A$ maps into
$K$ and since $\Sigma (E) = \Sigma (K) \oplus B^{**}$, any
extensions of the maps $\varphi_i^K$ to $E_{x_i}$ will satisfy
the requirements of the lemma.

Suppose now that the extension does split, so $B$ is in fact
locally free.  Fix maps $\varphi_i^K$ as above and choose maps
$\varphi_i^B \colon B_{x_i} \onto {\m{C}}$; set $\varphi_i :=
\varphi_i^K + \varphi_i^B $.  If $Hom(B,K) = 0$ then for any $A$
as in the last paragraph, the map $A \to B$ is either $0$ or an
isomorphism and $A \to K$ is then an inclusion or $0$
respectively, so the conclusion of the lemma is again satisfied.
If, on the other hand, there does exist a non-zero homomorphism
from $B$ into $K$, let $mb -c$ be the multiplicity of $B$ in
$\Sigma (K)$ for some integers $m,c$ with $0 \le c < b$.  By the
inductive hypothesis, for any non-zero map $B \to K$ the image of
$B$ is contained in kernels of at most $m-1$ of the maps
$\varphi_i^K$.  Suppose then that for each choice of maps
$\varphi_i^B \,\colon B_{x_i} \onto {\m{C}}$ there is a map $B
\to K$ such that {\it more\/} than $m-1$ of the maps $\varphi_i$
restricted to the image of $B$ are identically $0$; let there be
$m+k$ of them generically.  Such a map $B \to K$ is uniquely
determined since the non-zero difference between any two would
give a map $B \to K$ whose images in $K_{x_i}$ would be contained
in more than $m-1$ of the subspaces $ker\,\varphi_i^K$.  Thus one
obtains a map $\,\Psi \,\colon\, \bigoplus_{i=1}^n B_{x_i}^* \to
\Gamma(X,Hom (B,K))$ (which is clearly linear) such that, for
some $i_1,\ldots,i_{m+k}, $ the compositions $\varphi_i^K \circ
\Psi (\varphi^B_1,\ldots,\varphi^B_n) + \varphi_i^B , i =
i_1,\ldots,i_{m+k}$ are all identically $0$ as maps $B \to
{\m{C}}$.  Equivalently, $\Psi$ can be viewed as a linear map
$B\otimes \bigoplus_{i=i}^n B_{x_i}^* \to K$ such that, after
restriction to the points $\{x_i\}$, the composition
$End(B_{x_i}) = B_{x_i}\otimes B_{x_i}^* \to K \to {\m{C}}$ is
just $-trace,\, i= i_1,\ldots i_m$.  If, for some $i$ the map
$B\otimes B_{x_i}^* \to K$ annihilates every section restricting
to a trace-free endomorphism at $x_i$, then it is easy to see
that $B$ must have rank $1$.  Hence regardless of the rank of $B$
there is an inclusion $B\otimes \bigoplus_{j=1}^{m+k}
B_{x_{i_j}}^* \to K$, giving an embedding of $(m +k)b$ copies of
$B$ in $K$.  Since the multiplicity of $B$ in $K$ is $mb -c$, it
follows that $c = 0 = k$ and the multiplicity of $B$ in $E = K
\oplus B$ is $mb +1 = (m+1)b -(b-1)$.  Moreover, for generic
$\varphi_i^B$, the homomorphisms $\varphi_i$ have the required
property. \qed

\bigskip
Suppose that $E$ is as above and that $\Sigma (E) =
\bigoplus_i V_i \otimes A_i$ where $V_i$ is a $d_i$-dimensional
vector space and $A_i$ is a stable $a_i$-bundle with $\mu (A_i) =
\mu (E)$ and with $A_i \not\simeq A_j$ for $i \not= j$.  Pick any
$n$ points $x_1,\dots,x_n \in X\backslash B(E)$ and glue in the
bundle $E_1$ at each of these points in the generic fashion
described by Lemma~4.1 to obtain a bundle $\t E$ on the blowup
$\blowup X$ of $X$ at all $x_1,\dots,x_n$.  If $\sigma_i$
represents $-L_i$ let $\rho := \sum_i \sigma_i$ and for $\epsilon
>0$ take $\omega_{\epsilon} := \pi^*\omega + \epsilon\rho$.

Let $\t A \subset \t E$ have rank $a$, have torsion-free
quotient, and maximize $\mu (*,\omega_{\epsilon})$ over the
admissible subsheaves of $\t E$ by $\t A$ for all $\epsilon > 0$
sufficiently small.  Since $\mu (\t A,\omega_{\epsilon}) = \mu
(\pi_*\t A) + \epsilon\rho \!\cdot\!c_1(\t A)/a $, letting
$\epsilon \to 0$ shows that $\mu (\pi_*\t A) = \mu (E)$.  By the
discussion preceding Lemma~4.1, $\sigma_i \!\cdot\! c_1(\t A)$ is
either $0$ or $-1$ according respectively to whether or not
$(\pi_*\t A)^{**}$ is contained in the kernel of the linear form
defining the gluing of $E$ to $E_1$ at $x_i$.  By Lemma~4.1
itself, the former occurs at most $\kappa := \max_i \{d_i/a_i\}\;
(\le r) $ times.  Since $\mu (\t E,\omega_{\epsilon}) = \mu (E)
-\epsilon n/r$, it follows that $\mu (\t A,\omega_{\epsilon}) <
\mu (\t E,\omega_{\epsilon})$ if $n > r\kappa$, in which case $\t
E$ is stable with respect to $\omega_{\epsilon}$ for all
$\epsilon > 0$ sufficiently small.

\medskip
To summarise:

\bigskip
\noindent {\bf Proposition~4.2.} \quad {\sl Let $E$ be a
semi-stable $r$-bundle with $\Sigma (E) = \bigoplus_i V_i \otimes
A_i$ where $V_i$ is a $d_i$-dimensional vector space and $A_i$ is
a stable $a_i$-bundle with $\mu (A_i) = \mu (E)$ and with $A_i
\not\simeq A_j$ for $i \not= j$.  If $n > r[\max_i \{d_i/a_i\}]$
then for any choice of $n$ points in $X \backslash B(E)$ there is
a bundle $\t E$ on the blowup $\blowup X$ of $X$ at these points
such that $\t E$ restricts to ${\cal O}(1) \oplus \sum_1^{r-1}
{\cal O}$ on each component of the exceptional divisor, $(\pi_*\t
E)^{**} = E$, and $\t E$ is stable with respect to
$\omega_{\epsilon} := \pi^*\omega + \epsilon\sum_{i=1}^n \sigma_i
$ for all $\epsilon >0$ sufficiently small. \quad } \qed

\bigskip
\noindent{\bf Remarks:} \quad {\bf 1. }  Since
$\pi_*E(-1)$ is locally free, there is a map $\pi^*E\to \t E(-1)$
obtained by pulling-back and composing with the canonical map
$\pi^*\pi_*\t E(-1) \to \t E(-1)$.  Similarly, since $\pi_*(\t
E^*)$ is locally free, (after dualising) there is a map $\t E \to
\pi^*E$, and the composition $\pi^*E(1) \to \t E \to \pi^*E$ is
$s\,{\bf 1}_{ E}$ where $s \in \Gamma({\cal O}(-1))$ defines the
exceptional divisor.

\smallskip
\noindent{\bf 2. } If $r=2$ the entire procedure is
considerably simplified.  In particular, if $E_0$ is given by an
extension $0 \to L_1 \to E_0 \to L_2\otimes {\cal J}\to 0$ for
some line bundles $L_1,\,L_2$ with $deg(L_1)=deg(L_2)$ and some
sheaf of ideals ${\cal J}$ with ${\cal O}/{\cal J}$ supported at
a finite set, then if $L_1 \not\simeq L_2$ or the extension does
not split, it suffices to blow up at one point and take $\t E$ to
be a non-trivial extension $0 \to \pi^*L_1(1) \to \t E \to
\pi^*L_2\otimes {\cal J} \to 0$.  If $E_0 = L_1\oplus L_1$ then
it suffices to blow up at $3$ points and take a non-zero
extension $0 \to \pi^*L_1(1,-1,1) \to \t E \to \pi^*L_1(0,0,0)
\to 0$.

\smallskip
\noindent{\bf 3. }  Instead of gluing-in the bundle
$E_1 = {\cal O}(1) \oplus {\cal O}^{r-1}$ to construct $\t E$,
another natural choice is to glue in the bundle ${\cal
O}(1)\oplus {\cal O}(-1)  \oplus {\cal O}^{r-2}$ (assuming of
course that $r \ge 2$).  The new bundle now has the same first
Chern class (and determinant) as $E_0$, and charge one unit
greater (as opposed to $(r-1)/2r$ units greater).  An analysis
similar to that which was given above to prove Proposition~4.2
should be possible in this case, but the calculations are more
involved and have not as yet been carried out.

\bigskip

As mentioned in the introduction to this section, the gluing
construction also yields a simple way to by-pass the technical
difficulties of [D5] which can occur because of potential
singularity of moduli spaces of stable bundles.

For a stable $r$-bundle $E$ on the compact surface $(X,\omega)$
with Hermitian-Einstein connection $\nabla$, the cokernel of
$\nabla_+ \colon \Lambda^1(End_0\, E) \to \Lambda_+^2(End_0\, E)$
vanishes iff $H^2(X,End_0\, E)$ does, where $End_0$ denotes
trace-free endomorphisms. Thus $\nabla$ is a smooth point in the
moduli space of irreducible Hermitian-Einstein connections iff
$E$ is a smooth point in the moduli space of stable bundles. By
Serre duality, $H^2(X,End_0\, E) \simeq H^0(X,End_0\, E\otimes
K_X)^*$ where $K_X$ is the canonical bundle of $X$.  If $s \in
H^0(X,End_0\, E\otimes K_X)$ is a non-zero section, pick a point
$x_0$ at which $s$ is not zero, and blowup $X$ at $x_0$. Now take
$E_0 = E$, $E_1 = {\cal O}(-1)\oplus {\cal O}(1) \oplus {\cal
O}^{r-2}$ and construct a bundle $\t E = E_0 \#_{\rho} E_1$ on
$\t X$ as in \S 2.  Since $K_{\t X}\simeq (\pi^*K_X)(-1)$ ([BPV,
Theorem~I.9.1]),  Lemma~2.3(a) implies that $\pi_*(End_0\,\t
E\otimes K_{\t X}(-1))$ is locally free, and therefore it is
isomorphic to $End_0\,E\otimes K_X$ since it agrees with this
sheaf away from $x_0$.  The direct image of the sequence
$$
End_0\,\t E\otimes K_{\t X} 0 \to End_0\,\t E\otimes K_{\t X} \to
End_0\,\t E\otimes K_{\t X}(-1) \to End_0\,\t E\otimes K_{\t
X}\otimes{\cal O}_{L_0}(-1) \to 0
$$
thus gives $0 \to
\pi_*End_0\,\t E\otimes K_{\t X} \to End_0\, E \otimes K_X \to \m
C_{x_0}$, and it is straightforward to check that the under the
composition $H^0(X,End_0\,E\otimes K_X) \simeq H^0(\t X,End_0\,
\t E \otimes K_{\t X}(-1)) \to H^0(L_0,\t E \otimes K_{\t X}(-1))
\simeq \m C$ the section $s$ is not mapped to zero for generic
$\rho$.   For such $\rho$, it follows that the dimension of
$H^2(\t X,End_0\,\t E)$ is one less than that of
$H^2(X,End_0\,E)$ so after performing this operation at enough
points the bundle $\t E$ on $\t X$ will satisfy $H^2(\t
X,End_0\,\t E)=0$. By semi-continuity of cohomology the same will
be true for any bundle on $\t X$ sufficiently near $\t E$, and by
Proposition~3.4 the bundle $\t E$ will be stable with respect to
$\omega_{\alpha}$ for all suitable $\alpha$ sufficiently close to
$0$.

More generally, if $L$ is any line bundle on $X$  the same
technique shows that by blowing up at least
$h^2((End_0\,E)\otimes L)$ generic points in $X$, the generic
bundle $\t E= E \#_{\rho} E_1$ on the blowup satisfies
$H^2((End_0\,\t E)\otimes \pi^*L)=0$.

\smallskip
To summarise:

\bigskip
\noindent {\bf Proposition~4.3.} \quad {\sl If $E$ is an
$r$-bundle on $X$ and $L$ is a holomorphic line bundle, then for
any set $T$ of $n \ge h^2(X,(End_0\,E)\otimes L)$ points in
general position there is a blowup $\t X$ of $X$ centered at $T$
together with a bundle $\t E$ on $\t X$ restricting to ${\cal
O}(1)\oplus {\cal O}(-1)\oplus {\cal O}^{r-2}$ on each component
of the exceptional divisor satisfying $(\pi_*\t E)^{**}=E$ and
$H^2(\t X,(End_0\,\t E)\otimes \pi^*L)=0$. \quad \qed}

\bigskip
If $X'$ is a blowup of $X$ with exceptional divisor
$D'$, and $E'$ is a bundle on $X'$ which is topologically trivial
on $D'$, the cokernel of $H^1(X', End_0\, E') \to H^1(D', End_0\,
E')$ is the kernel of the epimorphism $H^2(X',(End_0\,E')(-D'))
\to H^2(X',End_0\,E')$.  If the former group vanishes, every
small deformation of $E'$ on $D'$ is induced by a small
deformation of $E'$ on $X'$.  Since every bundle on $\m P_1$ with
$0$ first Chern class has arbitrarily small deformations which
are holomorphically trivial, it follows in this case that $E'$
has arbitrarily small deformations which are holomorphic
pull-backs from $X$.

If $deg(K_X,\omega)<0$ and $E'$ is $\omega_{\alpha}$-stable for
sufficiently small suitable $\alpha$, $End_0\, E'\otimes
K_{X'}(D')$ has negative degree with respect to $\omega_{\alpha}$
and therefore $H^2(X',(End_0\,E')(-D'))=0$ by duality.  If
$deg(K_X,\omega)\ge 0$, applying Proposition~4.3 (blowing up
points of $X'\backslash D'$) yields a blowup $\t X'$ of $X'$ with
exceptional divisor $D$  together with a bundle $\t E'$ on $\t
X'$ such that $H^2(\t X',(End_0\,\t E')(-D'))=0$. From the
previous paragraph it follows that there are arbitrarily small
deformations of $\t E'$ which restrict to holomorphically trivial
bundles on $D'$.  The behaviour of the deformations on $D$ cannot
be controlled, but for sufficiently small perturbations of $\t
E'$, the bundles must be either trivial or split as  ${\cal
O}(1)\oplus{\cal O}(-1)\oplus {\cal O}^{r-2}$ on each component
by semi-continuity of cohomology.

\bigskip
\bigskip
\bigskip

\noindent {\bf 5.\quad Structure of moduli spaces II.}

\bigskip
Let $X$ be a compact complex surface, equipped with a
$\dbd$-closed positive $(1,1)$-form $\omega$. Consider a sequence
$\{E_i\}$ of stable holomorphic bundles of  fixed topological
type and degree, identified with a sequence $\{A_i\}$ of
Hermitian-Einstein connections on a fixed unitary bundle.  For
this sequence of connections, the $L^2$ norms of the curvatures
are uniformly bounded, so by Uhlenbeck's theorem [Sed], [U1]
there is a finite set $S \subset X$ and gauge transformations
such that the gauge-transformed sequence (also denoted $\{A_i\}$)
converges weakly in $L^2_{1,{\rm loc}}(X\backslash S)$ to a
connection $A$ defining a finite-action Hermitian-Einstein
connection.  Ellipticity of the Hermitian-Einstein equations
combined with Donaldson's argument in the proof of Corollary~23
of [D2] shows that a subsequence is converging weakly in
$L^p_{1,{\rm loc}}(X\backslash S)$ for any $p$ and by
bootstrapping and diagonalisation, that a subsequence converges
strongly in $C^k$ on compact subsets of $X\backslash S$ for any
$k$.  Since a Hermitian-Einstein connection can be twisted
locally by a Hermitian-Einstein connection on a trivial bundle so
that the new connection has $\lambda = 0$, it follows from the
Removability of Singularities theorem [U2] that the limit extends
across $S$ to define a new Hermitian-Einstein connection on a
bundle over $X$, and therefore a new semi-stable bundle $E$ with
$\Sigma(E) = E$.  The new holomorphic bundle $E$ has the same
rank and first Chern class as the bundles in the sequence, its
determinant is the limit of the determinants, but its charge is
at least one less for each point in $S$ where the curvature has
``bubbled".

Following [D5], this type of ``convergence" for sequences of
connections is referred to as {\it weak convergence\/} (on $X
\backslash S$), and sequences of stable holomorphic bundles of
the same degree and topological type converge weakly (on
$X\backslash S$) (with respect to $\omega$) if the  corresponding
irreducible Hermitian-Einstein connections converge  weakly.

In dealing with limits of sequences of stable bundles, arguments
are greatly simplified whenever it is known that a weak limit is
itself stable, rather than just semi-stable.  The stabilisation
construction given in the previous section is designed to meet
this type of need, and when combined with Lemma~2.2 of [B4]
(semi-continuity of cohomology on $H^0$ for weak limits), the
upshot is Lemma~5.1 below.  Recall that the notation $\Lambda\,F$
denotes $*(\omega \wedge F)$ and if $A$ is a connection on a
unitary bundle and $g$ is a complex automorphism (an {\it
intertwining operator}) on that bundle, $g \!\cdot \! A$ is the
connection with $(0,1)$ part $\db_{g\cdot A}= g\circ \db_A\circ
g^{- 1}$ and $(1,0)$ part $\d_{g\cdot A}=(g^*)^{-1}\circ
\d_A\circ g^*$.

\bigskip
\noindent{\bf Lemma~5.1.} \quad {\sl Let $\{E_i\}$ be a
sequence of stable bundles on $X$ corresponding to a sequence
$\{A_i\}$ of Hermitian-Einstein connections on a fixed
$U(r)$-bundle and let $S\subset X$ be a finite set such that the
sequence converges weakly on $X\backslash S$ to a
Hermitian-Einstein connection $A$ defining a quasi-stable bundle
$E$.  If $(\t X,\t E)$ is a stabilisation of $E$ with the none of
the blown-up points lying in $S$, then there is a subsequence
with stabilisations $\t E_{i_j}$ on $\t X$ converging weakly to
$\t E$ on $\t X \backslash \pi^{-1}(S)$. Here stability on $\t X$
is with respect to $\omega_{\epsilon}=\omega_{\epsilon\alpha_0}$
for $\alpha_0 := (1,\dots,1)$ and any $\epsilon>0$ sufficiently
small.}

\medskip
\noindent {\bf Proof:} \quad Let $T \subset X
\backslash B(E)$ be the finite set used to stabilize $E$.  Since
there is no bubbling of curvature near the points of $T$, a
sequence of integrable connections $\t A_i$ corresponding to
bundles $\t E_i$ with $\pi_*\t E_i = E_i$ can be found such that
$\t A_i$ agrees with $\pi^*A_i$ outside a fixed neighbourhood $\t
U$ of the exceptional divisor and which converge smoothly inside
this neighbourhood to a connection inducing $\t E$ there, so the
$\t A_i$ converge weakly in $L^p_{1,{\rm loc}}(\t X \backslash
\pi^{-1}(S))$ to a connection inducing $\t E$.

If $b_1(X)$ is even then by Proposition~3.5, all of the bundles
$\t E_i$ and $\t E$ will be stable with respect to the same
metric $\omega_{\epsilon}$ for sufficiently small $\epsilon>0$.
However, by construction of $\t E$ and $\t E_i$, a reflexive
subsheaf $\t A_i$ minimising $\nu_{\t E_i}(*,\omega_{\epsilon})$
for all sufficiently small $\epsilon$ and sufficiently large $i$
will satisfy the same type of splitting behaviour as that
described in the discussion preceding Proposition~4.2, and
therefore there exists $\epsilon_0>0$ such that $\t E_i$ is
$\omega_{\epsilon}$-stable for all $i$ and all
$\epsilon<\epsilon_0$, regardless of the parity of $b_1(X)$.

By weak compactness on $\t X$ there is a finite set $\t S \subset
\t X$ such that, after gauge transformations a subsequence of the
corresponding $\omega_{\epsilon}$-Hermitian-Einstein connections
$\t A'_i$ inducing $\t E_i$ converges weakly in $L^p_{1,{\rm
loc}}(\t X \backslash \t S)$ to a Hermitian-Einstein $\t A'$
connection on $\t X$ which, after removal of singularities,
defines a semi-stable bundle $\t E'$ there.  By Lemma~2.2 of
[B4], after rescaling if necessary the automorphisms $g_i$
intertwining $\t A'_i$ with $\t A_i$ give rise to a non-zero
holomorphic map $\t E \to \t E'$, but since the former is stable
and the two bundles have the same degree, it must be an
isomorphism.  This implies that the sequences $\{g_i^{}\},
\{g_i^{-1}\}$ are both uniformly bounded in $C^0(\t X)$, and by
Lemma~2.1 of that same reference it follows that $\t S =
\pi^{-1}(S)$ and that the two sequences bubble the same amount of
charge at each point of $S$.  \quad \qed

\bigskip

{}From the complex analytic viewpoint, Propositions 3.4 and 3.5
indicate that on a blowup $\t X$ of $X$ there exist ``stable"
moduli spaces ${\cal M}(\t X,E_{\rm
top},\omega_{\epsilon\alpha})$ of stable bundles as $\epsilon$
tends to $0$, at least when $b_1(X)$ is even.  Given a bundle $\t
E$ in one of these moduli spaces, it is natural to enquire about
the behaviour of the corresponding sequences of
Hermitian-Einstein connections and to determine the extent to
which the behaviour of this sequence reflects the complex
analytic picture of a bundle on $X$ glued to a bundle on a
neighbourhood of the exceptional divisor according the
prescription of \S 2.  These questions are at least partially
answered by Corollary~5.3 below where it is shown that the
corresponding connections converge off the blown-up set and
$B(\pi_*E)$ to the Hermitian-Einstein connection on
$\Sigma(\pi_*E)$;  the behaviour of the sequence near the
exceptional divisor is described in Proposition~5.4.

\bigskip
\noindent {\bf Proposition~5.2.}\quad {\sl Let $\blowup
X$ be a modification of $X$ with exceptional divisor $D \subset
\t X$, and let $\{\omega_i\} $ be a sequence of
$\bar\partial\partial$-closed positive $(1,1) $-forms on $\t X$
converging smoothly to $\pi^*\omega$.  Let $\t E$ be a
holomorphic $r$-bundle on $\t X$ and let $\{A_i\}$ be a sequence
of smooth hermitian connections on $\t E$ such that
\medskip
\itemitem{(i)\quad} $\norm{\Lambda_i F(A_i)
-\sqrt{-1}\lambda_i\,{\bf 1}}_{L^1(\t X, \omega_i)}\,\to \,0\,,$
for $\lambda_i = -2\pi\mu (\t E,\omega_i)/ Vol\,(\t X,\omega_i)$,
and
\itemitem{(ii)\quad} there is a finite set $S \subset X$, $p
>4$ and $k\ge 1$ such that $\pi_*A_i$ converges weakly in
$L^p_{k,{\rm loc}}(X\backslash (S \cup \pi (D)),\omega)$ to a
connection $A_{\infty}$ with $\norm{F(A_{\infty})}_{L^2(X)} <
\infty$.

\smallskip
\noindent Then

\medskip
\itemitem{(a)\quad} $\pi_*\t E$ is semi-stable;
\itemitem{(b)\quad} the quasi-stable bundle $E_{\infty}$ on $X$
defined by $A_{\infty}$ (after removal of singularities) is
isomorphic to $\Sigma(\pi_*\t E)$;
\itemitem{(c)\quad} after
suitable gauge transformations, a subsequence of $\{\pi_*A_i\}$
converges weakly in $L^p_{k}$ and strongly in $C^{k-1}$ on
compact subsets of $X\backslash \big(\pi(D) \cup B((\pi_*\t
E)^{**})\big)$.\par}

\medskip
\noindent {\bf Proof:}\quad If $L$ is an irreducible
component of $D$, the curvature form $f_i$ of the
Hermitian-Einstein connection on ${\cal O}(L)$ corresponding to
the metric $\omega_i$ satisfies $\norm{f_i}_{L^2(\t X,\omega_i)}
= 4\pi^2\,[1-\,{\rm deg}({\cal O}(L),\omega_i)^2/{\rm Vol}\,(\t
X,\omega_i)]$, which converges to $4\pi^2$.  Ellipticity of the
Hermitian-Einstein equations and the convergence of the sequence
$\{\omega_i\}$ to $\pi^*\omega$ implies that (a subsequence of)
$\{f_i\}$ converges weakly in $L^2(X)$ and smoothly on  compact
subsets of $X\backslash\pi(L)$ to a finite action solution of
the Hermitian-Einstein equations which, by removable
singularities,  extends to $X$ to define a holomorphic line
bundle of degree $0$  there.  By Lemma~2.2 of [B4] and Hartogs'
theorem, this line bundle has a  non-zero holomorphic section,
but since the degree of the bundle is $0$, this section is
covariantly constant and therefore defines a global
trivialisation---all of the curvature has concentrated along $L$.
It follows from this that if $\t E$ is twisted by a suitable line
bundle trivial off $D$ and the sequence $\{A_i\}$ is
correspondingly twisted by the sequence of Hermitian-Einstein
connections on these line bundles, then the hypotheses of the
proposition remain true for the new sequence and it can be
supposed from now on that $\pi_*\t E$ is locally free.

By Lemma~2.2 of [B4] and Hartogs' theorem again there is a
non-zero holomorphic map $E_{\infty} \to \pi_*\t E$, so if either
bundle is stable this map is an isomorphism.  In particular, this
occurs if $\t E$ is a line bundle.

\medskip
\noindent (a) \quad If $A \subset \pi_*\t E$ has rank
$a$ and torsion-free quotient  then $ \nu_{\pi_*\t E}(A) =
\nu_{\t E}(\pi^*A,\omega_i) -{1\over r}(\omega_i
-\pi^*\omega)\!\cdot\! {\sl det}\,(\t E\otimes \pi^*A^*)$.  Since
$\nu_{\t E}(\pi^*A,\omega_i)\ge \nu_{\t E}(\h A,\omega_i)$ where
$\h A$ is the maximal normal extension of $\pi^*A$ in $\t E$ and
the latter term has non-negative limit by hypothesis {\sl (i)\/}
and (1.7), the convergence of $\omega_i$ to $\pi^*\omega$ implies
$ \nu_{\pi_*\t E}(A) \ge 0$.

\medskip
\noindent (b) \quad If either $E_{\infty}$ or $\pi_*\t
E$ is stable, this is proved already by the remarks above. If
$\pi_*\t E$ is not stable, there exists a subsheaf $A \subset
\pi_*\t E$ with torsion free quotient $C$ such that $\nu_{\pi_*
\t E}(A) = 0$.  Pulling back to $\t X$ and taking the maximal
normal extension gives a subsheaf $\h A \subset \t E $ with
torsion-free quotient $\h C$ such that $\pi_*\h A = A$ and $C
\subset \pi_*\h C$.  Desingularise the sequence $0 \to \h A \to
\t E \to \h C \to 0$ as in \S 3 of [B2] to obtain a modification
$\t X' \oa {\pi'} \t X \oa \pi X$ and a sequence of bundles $0
\to \h A' \to \pi'^*\t E \to \h C' \to 0$ on $\t X'$ with
$\pi'_*\h A' = \h A$ and $\h C \subset \pi'_*\h C'$.

Choose metrics $\omega_i' = \pi'^* \omega_i + \rho_{\alpha_i} $
on $\t X'$ with $\rho_{\alpha_i} $ converging smoothly to $0$,
and make this convergence sufficiently fast so that
$\omega_i\wedge'F(\pi'^*A_i) -\sqrt{-1}\lambda'_i\,{\bf
1}\omega_i'^2$ converges to $0$ in $L^1(\t X',\omega'_i)$; this
is possible since the connections $A _i$ are smooth (cf.~the
final remark in \S 2 of [B2], p.631).  From (1.7) and {\sl (i)\/}
it follows that if $\beta_i$ is the second fundamental form of
$\h A'$ in $\pi'^*\t E$ for the hermitian connection $\pi'^*A_i$,
then $\beta_i$ is converging to $0$ in $L^2(\t X',\omega'_i)$.
This implies that the hypotheses of the proposition are satisfied
by the induced connections on $\h A'$ and $\h C'$.  If
$A_{\infty},\, C_{\infty} $ are the holomorphic bundles on $X$
defined by the limiting (Hermitian-Einstein) connections on $\h
A',\,\h C'$, then since $\beta_i \to 0 $ it follows from the
construction that $E_{\infty} = A_{\infty} \oplus C_{\infty}$.
By induction on rank, $A_{\infty} = \Sigma\big((\pi\pi')_*\h
A'\big)$ and $C_{\infty} = \Sigma\big((\pi\pi')_*\h C'\big)$,
giving $ E_{\infty} \,=\, A_{\infty} \oplus C_{\infty} \,=\,
\Sigma\big((\pi\pi')_*\h A'\big) \oplus \Sigma\big((\pi\pi')_*\h
C'\big) \,=\, \Sigma(A) \oplus \Sigma(C) \,=\, \Sigma(\pi_*\t
E)$, proving (b) in general.

\medskip
\noindent (c)\quad By passing to a subsequence it can be
supposed that the sequence $\{g_i\}$ of intertwining operators
from which the holomorphic map $\pi_*\t E \to E_{\infty}$ is
constructed is converging strongly in $C^k$ on any compact subset
of  $X \backslash (S \cup \pi(D))$. If $\pi_*\t E$ is stable then
this map  is an isomorphism and the sequence $\{g_i^{-1}\}$ must
also be  bounded in $C^k$ on compact subsets of $X\backslash (S
\cup\pi(D))$.  Then  Lemma~2.1 of [B3] implies that curvature can
only concentrate on $\pi(D)$ itself,  which proves (c) in this
case.  The general case now follows by  induction on rank using
the proof of (b) when $\pi_*\t E$ is not  stable, since any point
of $\t X$ in the center of the  modification $\t X' \to \t X$ is
mapped into $B(\pi_*\t E)$ by $\pi$. \quad\qed

\bigskip

\noindent {\bf Corollary~5.3.}\quad {\sl If $\t E$ on $\t X$ is
stable with respect to $\omega_{\alpha_i}$ for some suitable
$\alpha_i \in \m R^n_+$ converging to $0$, then after suitable
gauge transformations, a subsequence of the corresponding
sequence $\{\nabla_i\}$ of Hermitian-Einstein connections
converges smoothly on compact subsets of $X \backslash  (\pi(D)
\cup B((\pi_*\t E)^{**}))$ to a Hermitian-Einstein connection
inducing $\Sigma(\pi_*\t E)$.}

\medskip
\noindent {\bf Proof:} \quad  On any compact subset of
$X\backslash \pi(D)$ the $L^2$ norm of $F(\nabla_i)$ with respect
to $\omega$ is uniformly bounded.  The weak compactness arguments
of [Sed], [U1] still apply in this setting to obtain a
subsequence of $\{\nabla_i\}$ (after gauge transformations)
converging weakly in $L^2_{1}$ on compact subsets of $X\backslash
(S\cup \pi(D))$ for some finite set $S \subset X$.  Ellipticity
of the Hermitian-Einstein equations together with boot-strapping
and diagonalisation then give a subsequence converging smoothly
on compact subsets of this complement, and the conclusion then
follows from the proposition. \quad \qed

\bigskip
Corollary~5.3 describes the behaviour of sequences of
Hermitian-Einstein connections away from the exceptional divisor
$D$.  To describe the behaviour {\it near\/} $D$ consider the
case that $\t X$ is the blowup of $X$ at a single point $x_0$ and
fix a metric on $\t X$ of the form $\omega_{\epsilon} =
\pi^*\omega + \epsilon \sigma$ where $\sigma$ restricts to the
Fubini-Study metric on $L_0 := \pi^{-1}(x_0)$.  Then as $\epsilon
\to 0$ the corresponding sequence of Hermitian-Einstein
connections on the line bundle ${\cal O}(1)$ over $\t X$
converges smoothly on compact subsets of $\t X \backslash L_0$ to the
trivial flat connection.

Near $L_0$ on the other hand, restrict attention to a small
neighbourhood of $x_0$ isomorphic to $B_{r_0}$ and choose
holomorphic coordinates $(z^0,z^1)$ in that neighbourhood so that
$\omega(x_0)$ is standard in these coordinates. If $\omega$ is
the standard flat metric $\delta :=(i/2)\ddb |z|^2$ near $x_0$
then under the dilation $z \mapsto \sqrt{\epsilon}z$ the metric
$\delta_{\epsilon}$ pulls-back to $\epsilon\delta_1$; in general,
if $\omega$ is arbitrary the pull-back of $\omega_{\epsilon}$
differs from $\epsilon\delta_1$ by a term of order
$\epsilon^{3/2}$ ($\epsilon^2$ if $d\omega(x_0) = 0$ and the
coordinates are appropriately chosen) as $\epsilon\to 0$.

{}From the local description of bundles on a blowup given in
Proposition~2.5 it is easy to see that the pull-backs of a
natural transition function for $\t E$ under the maps
$\lambda_{\epsilon}$ given by $ B(x_0,r_0/\sqrt{\epsilon}) \owns
z \mapsto \sqrt{\epsilon}z\in B(x_0,r_0) $ converge smoothly to a
transition function for a direct sum of line bundles on $\t{\m
C}{}^2$.

\smallskip
The following result can be proved using Lemma~5.1,
Corollary~5.3 and the  stabilisation procedure of \S4:

\bigskip
\noindent {\bf Proposition~5.4.} \quad {\sl Let $\blowup
X$ be the blowup of $X$ at $x_0$ with $L := \pi^{-1}(x_0)$, and
let $\t E$ be a bundle with $\t E\!\!\mid_L = \oplus_i {\cal
O}(a_i)$.  If $\t E$ admits an
$\omega_{\epsilon}$-Hermitian-Einstein connection
$\nabla_{\epsilon}$ for all $\epsilon > 0$ sufficiently small,
then as $\epsilon \to 0$ the pull-back of $\nabla_{\epsilon}$ to
$\t B(x_0,r_0/\sqrt{\epsilon})$ converges smoothly on compact
subsets of $\m C^2$ after suitable gauge transformations to the
direct sum of the standard anti-self-dual connections on ${\cal
O}(a_i)$, where anti-self-duality is with respect to the standard
metric $\delta_1$ on $\t{\m C}{}^2$ in suitable holomorphic
coordinates. \qed}

\medskip
\noindent (Cf.~[B3] for an interpretation of connections
of this form in terms of intantons on $\m{CP}_2$.)

\bigskip
\bigskip
\bigskip

\noindent{\bf 6. \quad Compactification of moduli spaces.}

\bigskip
\medskip

Let $X$ be a compact complex surface equipped with a positive
$\dbd$-closed $(1,1)$-form $\omega$, and let $\{E_i\}$ be a
sequence of stable bundles on $X$ converging weakly to a
quasi-stable bundle $E$. Blow-up $X$ along $S$ to $\blowup X$ and
pull back the sequence $\{E_i\}$. Assume temporarily that every
bundle in the pulled-back sequence is stable with respect to the
same metric $\omega_{\alpha}$ for some suitable $\alpha$
sufficiently small---this will be true if either $b_1(X)$ is even
(by Proposition~3.5) or if $E$ is stable. (If $\t A_i\subset
\pi^*E_i$ maximises $\mu(\bullet,\omega_{\epsilon\alpha})$ for
all $\epsilon>0$ sufficiently small, then $\pi_*\t A_i$ is
semi-stable.  If $\nu_{\pi^*E_i}(\t
A_i,\omega_{\epsilon_i\alpha})=0$ for $\epsilon_i \to 0$, then as
in \S3, the ranks, first Chern classes, charges, and degrees of
$\pi_*\t A_i$ are all uniformly bounded, and the same holds for
any stable $A_i' \subset \pi_*\t A_i$ of the same degree.
Passing to a subsequence, taking a weak limit and applying
semi-continuity of cohomology then yields a semi-stable bundle
$A$ with $rank(A)<rank(E)$, $\mu(A)=\mu(E)$ and a non-zero
holomorphic map $A\to E$, contradicting stability of $E$.)

By weak compactness on $\t X$ the  new sequence converges weakly
off some finite set $\t S \subset \t  X$ to a semi-stable bundle
$\t E$ and by semi-continuity of cohomology there is a  non-zero
holomorphic map $(\pi_*\t E)^{**} \to E$.  If either is stable,
then this map is an isomorphism and it follows from Lemma~2.1 of
[B4] that $\t S \subset \pi^{-1}S$;  moreover $C(\t E) \ge C(E)$
by Proposition~2.9.

This argument suggests that there is another compactification of
moduli spaces tied more closely to the complex analysis, distinct
from the more usual ones ([Gie], [M]): if $\blowup X$ is a blowup
of $X$ at $x_0$, sequences of stable bundles $\t E_i$ on $\t X$
are easily constructed such that each is trivial on the
exceptional line and such that the sequence converges smoothly on
$\t X$ to a bundle $\t E$ which is stable but which is now
non-trivial on the exceptional line.  Thus the theorems of
Uhlenbeck [U1], [U2] on compactness and removability of
singularities in the space of Yang-Mills connections are
reinterpreted as the well-known phenomenon of jumping of
holomorphic structures.  It is natural to conjecture that this is
the only type of catastrophe which can occur, and that, given a
degenerating sequence of Hermitian-Einstein connections on a
bundle over $X$, after a finite number blowups and pull-backs, a
strongly convergent subsequence can be found.  Unfortunately,
there is no guarantee that after blowing up and  pulling back,
the new sequence of Hermitian-Einstein connections  will not
bubble precisely the same amount of charge as the original  and
that iterating the process will eventually lead to a convergent
subsequence. Indeed, by a process of diagonalisation, sequences
of connections can be constructed for which this blowing-up and
pulling-back procedure will not terminate after finitely many
steps.  Despite this, provided that some flexibility is
introduced  into the blowing-up process, a form of compactness
{\it does\/}  hold, as is indicated by the following result
proved in [B4]:

\bigskip
\noindent{\bf Theorem~6.1.} \quad {\sl  Let $X$ be a
compact complex surface equipped with a positive $\dbd$-closed
$(1,1)$-form $\omega$.  Let $\{A_i\}$ be a sequence of
Hermitian-Einstein connections on a fixed unitary bundle $E_{\rm
top}$ of rank $r$ over $X$ such that the corresponding
holomorphic bundles $E_i$ are stable and are of uniformly bounded
degree.  Suppose that $E_i$ converges weakly to $E$ off $S
\subset X$.  Then for some subsequence $\{E_{i_j}\}\subset
\{E_i\}$:

\medskip
\item{1.\quad} There is a sequence of blowups $\t
X_{i_j}$ of $X$  converging to a blowup $\blowup X$ of $X$ with
exceptional divisor $\pi^{-1}(S)$ and with $c_2(\t X)\ le
c_2(X)+2C(E_{\rm top})-2C(E)-1$;
\item{2.\quad} There are complex
automorphisms $g_{i_j}$ of $\pi_{i_j}^*E^{}_{\rm top}$ such that
$|g_{i_j}^{}|+|g_{i_j}^{-1}|$ is uniformly bounded on  compact
subsets of $\t X\backslash \pi^{-1}(S)$ and $\{g_{i_j}\}$ is
converging uniformly in $C^3$ on  such subsets;
\item{3.\quad}
$\{g_{i_j}^{}\cdot (\pi_{i_j}^*A_{i_j}^{})\}$ converges uniformly
in $C^2$ to a smooth connection on $\pi^*E_{\rm top}$ over $\t X$
which defines a holomorphic bundle $\t E$ such that $(\pi_*\t
E)^{**}=E$;
\item{4.\quad} If $E$ is stable then for any suitable
$\alpha$ with $|\alpha|$ sufficiently small, the connections
$g_{i_j}^{}\cdot (\pi_{i_j}^*A_{i_j}^{})$ can be taken to be
$(\pi_{i_j}^*\omega+\rho_{\alpha})$-Hermitian-Einstein, and $\t
E$ is $\omega_{\alpha}$-stable. \par

\medskip
\noindent  If $E$ is not stable, but $b_1(X)$ is even
and $rank(E_{\rm top})=2$ there is sequence of blowups
$\{X_{i_j}\}$ consisting of at most $2C(E_{\rm top})-1$
individual blowups converging to a blowup $\t X$ such that, for
some suitable $\alpha$, $\pi_i^*E_i^{}$ is
$\omega_{\epsilon\alpha}$-stable for all $\epsilon \in (0,1]$ and
the subsequence $\{\pi_{i_j}^*E_{i_j}^{}\}$ converges strongly to
a bundle $\t E$ on $\t X$, stable with  respect to
$\omega_{\epsilon\alpha}$-stable for all $\epsilon\in (0,1]$. }

\medskip

\noindent The convergence of a sequence of blowups should be
interpreted as the convergence of a sequence of integrable
complex structures on the  same underlying smooth manifold $X\#n
\b{\m P}_2$ endowed with a fixed Riemannian metric.  In the
current setting, these complex structures are isomorphic on the
complement of  a fixed open set with strictly pseudoconvex
boundary.

A priori, there is no reason why the last statement of the
theorem cannot be extended to bundles of arbitrary rank, but as
yet a proof is lacking.  The complications arise, as always, from
the presence of bundles which are semi-stable but not stable.

\bigskip

Theorem~6.1 suggests a variety of different compactifications for
moduli of stable bundles---for example,  adding torsion-free
semi-stable sheaves which are direct images of stable bundles on
blowups is an obvious candidate.  One other which has some
interesting properties is described below.

\medskip

Let $E_{\rm top}$ be a unitary $r$-bundle over $(X,\omega)$ and
let ${\cal M}(X,E_{\rm top})$ denote the space of isomorphism
classes of quasi-stable holomorphic structures on $E_{\rm top}$.
Consider the set of pairs $(\t X, \t E)$ where

\smallskip
{\sl
\itemitem{1.\quad} $\t X$ is a blowup of $X$;
\itemitem{2.\quad} $\t E$ is a holomorphic bundle on $\t X$
topologically isomorphic to $\pi^*E_{\rm top}$ such that $\pi_*E$
is semi-stable;
\itemitem{3.\quad} $\t E$ is
$\omega_{\alpha}$-quasi-stable for all suitable $\alpha$ in an
open set of such;
\itemitem{4.\quad} If $b_1(X)$ is odd, $\deg(\pi_*\t E,\omega)=0$. \par}

\smallskip
\noindent Note that the third condition implies that
if $\t E$ is not stable, then it is a direct sum of stable
bundles each of which is topologically trivial on the exceptional
divisor.  Note also that the requirement that $\pi_*\t E$ be
semi-stable implies that $\t E$ is
$\omega_{\epsilon\alpha}$-(quasi-)stable for all $\epsilon \in
(0,1]$ if $\t E$ is $\omega_{\alpha}$-quasi-stable.

If $b_1(X)$ is odd, it follows from the proof of Proposition~2 in
[B2] that for any $t \in \m R$ there is a line bundle $L_t$ on $X$
with $c_1(L_t)=0$ and $\deg(L_t,\omega)=t$, so the fourth condition
simply prevents this type of non-compactness.

\medskip
\noindent  On the set of pairs $(\t X,\t E)$ satisfying
the above conditions, define an equivalence relation $\sim$ by
setting $(\t X_1,\t E_2) \sim (\t X_2,\t E_2)$ iff there is a
joint blowup $\t X_{12}$ such that  $\pi_2^*\t E_1\simeq
\pi_1^*\t E_2$ on $\t X_{12}$, and let  $\b{\cal M}(X,E_{\rm
top})$ denote the set of equivalence classes. A topology is
defined on $\b{\cal M}=\b{\cal M}(X,E_{\rm top})$ by defining
$\{[(\t X_i,\t E_i)]\}\subset \b{\cal M}$ to converge to $[(\t
X,\t E)]$ iff $[(\t X_i,\t E_i)]$ can be represented by a
sequence of blowups $\t X_i$ converging to $\t X$ with
$\omega_{\alpha}$-quasi-stable bundles $\t E_i$ on $\t X_i$
converging strongly to $\t E$ on $\t X$.

If $(\t X,\t E)\in \b{\cal M}$ with $\t E$ quasi-stable with
respect to $\omega_{\alpha'}$ for all $\alpha'$ in a
neighbourhood of $\alpha$ and $A_0$ is the
$\omega_{\alpha}$-Hermitian-Einstein connection on $E_{\rm top}$
inducing $\t E$,  it is easily verified that the image in
$\b{\cal M}$ of a set of the form
$U(A_0,\epsilon_1,\epsilon_2):=\{A_0+a\}$ where $A_0+a$ is
$\omega_{\alpha}$-Hermitian-Einstein,
$||a||_{C^1(\omega_{\alpha})} < \epsilon_1$ with $E(A_0+a)$
$\omega_{\alpha'}$-quasi-stable for $|\alpha'-\alpha|<\epsilon_2$
is open in $\b{\cal M}$ and that every open set containing $(\t
X,\t E)$ contains such a open neighbourhood.  Hence the
collection of such sets forms a base for the topology on $\b{\cal
M}$, from which it is clear that this topology is second
countable.

\smallskip
Theorem~6.1 implies that at least under certain
circumstances, every sequence in ${\cal M}$ has a subsequence
converging in $\b{\cal M}$.  To attempt to use this construct a
compactification of ${\cal M}$ clearly requires a bound on the
number of blowups required to represent classes in $\b{\cal M}$,
a bound which is provided by Proposition~6.3 below:

\bigskip
\noindent {\bf Lemma~6.2.} \quad {\sl For any $r_0, C_0$
there exists $n = n(r_0,C_0,\omega)$ such that any semi-stable
bundle $E$ of rank $\le r_0$ and charge $\le C_0$ satisfies
$h^2(X,End_0\,E) \le n$.}

\medskip
\noindent{\bf Proof:} \quad  If not, there exists a
sequence $\{E_i\}$ of semi-stable bundles of bounded rank and
charge such that $h^2(X,End_0\,E_i)\to \infty$.  It can be
assumed without loss of generality that each $E_i$ is
quasi-stable since it follows easily by induction on rank that
$h^2(End_0\,E) \le h^2(End_0\,\Sigma(E))$ for any semi-stable
$E$. The bundles $End_0\,E_i$ have uniformly bounded ranks ($\le
r_0^2-1$) and charges ($\le 2r_0C_0$) with trivial determinants,
and each is also quasi-stable.  By passing to a subsequence if
necessary, it can be assumed that $End_0\,E_i$ is converging
weakly to a quasi-stable bundle $B$ off some finite set $S
\subset X$.  Passing to another subsequence and using
Theorem~6.1, there is a sequence of blowups $\{\t X_i\}$
converging to a blowup $\t X$ such that a sequence of smooth
integrable connections representing $\pi_i^*End_0\,E_i^{}$
converges strongly in $C^2$ to define a bundle $\t B$ on $\t X$
with $(\pi_*\t B)^{**}=B$.   Then $h^2(\t
X_i,\pi_i^*End_0\,E_i^{})=h^0(\t
X_i,(\pi_i^*End_0\,E_i^{})\otimes K_{\t X_i})=
h^0(X,(End_0\,E_i)\otimes K_X)=h^2(X,End_0\,E_i) \to \infty$,
contradicting the fact that $h^2(\t X_i,End_0\,E_i) \le h^2(\t
X,\t B)$ for $i$ sufficiently large, by standard semi-continuity
of cohomology. \quad \qed

\bigskip

\noindent{\bf Proposition~6.3.} {\sl There is an integer
$N=N(X,\omega,E_{\rm top})$ such that every class in $\b{\cal M}$
can be represented by a bundle on a blowup consisting of at most
$N$ individual blowups.  If $\deg(K_X,\omega)<0$ then $N \le
2C(E_{\rm top})-1$.}

\medskip
\noindent{\bf Proof:}  Suppose $(\t X,\t E)\in \b{\cal
M}$ and $\t E$ is $\omega_{\alpha}$-quasi-stable for some
suitable $\alpha$, and let $\t S \subset \t X$ be the exceptional
divisor.  Pick a set of $r^2+1$ points $T_0 \subset X\backslash
\pi(supp(\rho_{\alpha}))$ and stabilise $(\pi_*\t E)^{**} =: E_0$
to a bundle $E_0'$ on a blowup $X'$ of $X$ centered at $T_0$.
Then for $\beta=\epsilon_0(1,1,\dots,1)$, it follows easily that
the corresponding ``stabilisation" $\t E'$ of $\t E$ on $\t X'$
is $\omega_{\epsilon\alpha,\delta\beta}$-stable for all
$\epsilon,\delta>0$ sufficiently small.

Since $\pi'_*(End_0\,E'_0(-1))=End_0\,E_0$, it follows that
$h^2(X',End_0\, E'_0)=h^2(X,End_0\,E_0)$.  Consequently, $h^2(\t
X',End_0\,\t E'\otimes {\cal O}(-\t S))=h^0(\t X',End_0\,\t
E'\otimes \pi^*K_{X'}\otimes {\cal O}(2\t S)) \le h^0(X',End_0\,
E_0'\otimes K_{X'}) = h^2(X,End_0\,E_0) \le n$, where $n$ is a
uniform bound as specified by Lemma~6.2.  If $\deg(K_X,\omega)<0$
then $H^2(X,End_0\, E)=0$ for any semi-stable $E$, so take $n=0$
in this case.

Pick a finite set $T_1 \subset \t X' \backslash
supp(\rho_{\alpha}+\rho_{\beta})$ consisting of at most $n$
points, blow up $X$ at these points to $X_1$ and construct a
bundle $E_1$ on $X_1$ such that $(\pi_{1*}E_1)^{**}=E_0$ with
$E_1$ restricting to ${\cal O}(-1)\oplus {\cal O}(1)\oplus {\cal
O}^{r-2}$ on every component of the exceptional divisor and with
$H^2(X_1,End_0\,E_1)=0$. There are corresponding bundles $ E'_1$
on the blowup of $X'$ along $T_1$ with $(\pi'_* E_1')^{**}= E_1$
and $\t E'_1$ on the blowup $\t X'_1$ of $\t X'$ along $T_1$
satisfying $(\pi_*\t E'_1)^{**}=E'_1$.  Moreover, if $\epsilon_0$
is fixed and sufficiently small with
$\alpha_0;=\epsilon_0\alpha$ and $\beta_0:=\epsilon_0\beta$, for
any $\gamma$ with $|\gamma|$ sufficiently small it follows  $\t
E_1'$ is $\omega_{\alpha_0,\beta_0,\gamma}$-stable and $ E'_1$ is
$\omega_{\beta_0,\gamma}$-stable.

By construction, $H^2(\t X_1',End_0\,\t E_1'\otimes {\cal O}(-\t
S))=0$, so the bundle $\t E_1'$ has arbitrarily small
deformations which are holomorphically trivial on $\t S$.  Hence
there is a sequence $\{E_{1,i}'\}$ of bundles on $X_1'$ such that
$\pi^* E_{1,i}'$ converges smoothly to $\t E_1'$.  Since this
bundle is $\omega_{\alpha_0,\beta_0,\delta\gamma}$-stable for all
$\delta \in (0,1]$, by passing to a subsequence it can be
supposed that the same is true of $\pi^*E_{1,i}'$, and therefore
$E_{1,i}'$ must be $\omega_{\beta_0,\delta\gamma}$-stable for
$\delta\in (0,1]$.

By passing to another subsequence if necessary, it can be assumed
that $\{E_{1,i}'\}$ is converging weakly to some
$\omega_{\beta_0,\delta\gamma}$-quasi-stable bundle $E_1'$ off
some finite set $T_2 \subset X_1'$, and by semi-continuity of
cohomology and $\omega_{\beta_0}$-stability of $E'_0$ it follows
that $(\pi_{1*}E_1')^{**}=E_0'$.   Since the sequence
$\{\pi^*E_{1,i}'\}$  converges smoothly to $\t E_1'$, it follows
from Lemma~2.1 of [B4] that $T_2 \cap \pi_1^{-1}(T_1)=\emptyset
=T_2 \cap\pi'^{-1}(T_0)$.

By construction, $C(E_{1,i}') \le C(E_{\rm top})+(r^2+1)/2r+n$,
and $C( E_1')$ must be at least $(r^2+1)/2r$ since $c_1(E_1')$
restricts to $1$ on each component of $\pi'^{-1}(T_0)$.  By
Theorem~6.1 therefore, after passing to another subsequence if
necessary there is a sequence of blowups $X'_{2,i}$ of  $X_1'$
converging to a blowup $X_2'$ centered at $T_2$ such that, if
$\pi_{2,i} \colon X_{2,i}' \to  X_1'$ is the blowing-down map,
$\pi_{2,i}^*E_{1,i}'$ converges to a bundle $E_2'$ on $X_2'$
satisfying $(\pi_{2*}E_2')^{**}= E_1'$,  with $c_2(X_2')\le c_2(
X_1') + 2(C(E_{\rm top})+n)-1$.  The convergence on this sequence
of blowups is with respect to a metric of the form
$\omega_{\beta_0,\gamma,\xi}$ for $|\xi|$ small;  fix one such
generic $\xi$.

A subsequence of the sequence of joint blowups $\t X_{2,i}'$ of
$X_{2,i}'$ and $\t X$ has a subsequence converging to a joint
blowup $\t X_2'$ of $X_2'$ and $\t X$.  Abusing notation
slightly,  the pull-back $\pi^*E_2'$ of $E_2'$ to $\t X_2'$ is
$\omega_{\epsilon\alpha_0,\beta_0,\gamma,\xi}$-stable for
$\epsilon>0$ sufficiently small, so by semi-continuity of
cohomology  the pull-backs of the bundles $E_{1,i}'$ must
converge to  $\pi^*E_2'$ with respect to
$\omega_{\epsilon\alpha,\beta_0,\gamma,\xi}$.  On the other hand,
the pull-back $\pi_2^*\t E_1'$ of $\t E_1'$ to $\t X_2'$ is
$\omega_{\alpha_0,\beta_0,\gamma,\delta\xi}$-stable for
$\delta>0$ sufficiently small, so with respect to this metric the
pull-backs of $E_{1,i}'$ converge to $\pi_2^*\t E_1'$.  By
semi-continuity of cohomology there is a non-zero holomorphic map
$\pi^*E_2'\to \pi_2^*\t E_1'$, giving a non-zero map
$\pi_{2*}E_2'=\pi_{2*}\pi_*\pi^*E_2' \to \pi_{2*}\pi_*\pi_2^*\t
E_1' =\pi_*\pi_{2*}\pi_2^*\t E_1'=\pi_*\t E_1'$.    On taking
double-duals, this gives a non-zero holomorphic map $E_1'\to
E_1'$;  since $E_1'$ is $\omega_{\beta_0,\gamma}$-stable, this
map is an isomorphism, and therefore, since
$\det(\pi^*E_2')=\det(\pi_2^*\t E_1')$, it follows that the map
$\pi^*E_2'\to \pi_2^*\t E_1'$ must be an isomorphism.

Let $X_3$ be the blowup of $X$ obtained by blowing down the
components of the exceptional divisor in $X_2'$ over $T_0$ and
$T_1$, so $c_2(X_3) =c_2(X_2')-n-(r^2+1) \le c_2(X)+2C(E_{\rm
top})+n-1$, and let $E_3 := (\pi_{1*}\pi'_*E_2')^{**}$.  Taking
double-duals of the direct image under $\pi'\pi_1$  of the
isomorphism $\pi^*E_2'\simeq \pi_2^*\t E_1'$ gives
$\pi^*E_3=\pi_2^*\t E$, so $(X_3,E_3)\sim (\t X,\t E)$ in
$\b{\cal M}$. Note that $\pi_2^*\t E$ is
$\omega_{\alpha,\delta\xi}$-quasi-stable for sufficiently small
$\delta>0$, which implies that $\pi^*E_3$ has the same property;
this in turn implies that $E_3$ is
$\omega_{\delta\xi}$-quasi-stable for small enough $\delta$.
Since $\xi$ is assumed to be generic, if $E_3$ is not actually
stable, each stable summand must be a topological pull-back from
$X$, so $E_3$ is $\omega_{\delta\xi'}$-quasi-stable for all
$\xi'$ near $\xi$. \quad \qed

\bigskip

\noindent{\bf Theorem~6.4.}
\medskip
{\sl

\item{1.~~} ${\cal M} \subset \b{\cal M}$ is an open set, dense
if $deg(K_X,\omega)<0$;
\item{2.~~} If $b_1(X)$ is even and
$rank(E_{\rm top})=2$ then $\b{\cal M}$ is compact;
\item{3.~~}
If every semi-stable bundle $E$ on $X$ satisfying
$rank(E)=rank(E_{\rm top})$, $c_1(E)=c_1(E_{\rm top})$ and $C(E)
\le C(E_{\rm top})$ is stable then $\b{\cal M}$ is a compact
Hausdorff space, smooth if $deg(K_X,\omega)<0$.\par}

\medskip
\noindent Note that if $b_1(X)$ is even and $c_1(E_{\rm
top})\in H^2(X,\m Z_r)$  is not a torsion class (where
$r=rank(E_{\rm top})$), the hypothesis in the third statement
holds for  generic $\omega$.  Even if $c_1(E_{\rm top})=0$, by
fixing a base point $x_0\in X$ and blowing up at this point, the
$\m P_{r-1}$-bundle over ${\cal M}_{stable} \subset {\cal M}$
with fibre $\m P(E_{x_0})$ at $E\in {\cal M}_{stable}$  is
isomorphic to the space of $\omega_{\epsilon}$-stable bundles on
$\t X$ restricting to ${\cal O}(-1)\oplus {\cal O}^{r-1}$ on
$\pi^{-1}(x_0)$ with direct image topologically isomorphic to
$E_{\rm top}$.

If $b_1(X)$ is odd, the hypothesis in the third statement is
quite a strong restriction since it will only be satisfied if
$c_2(E_{\rm top})<0$, for if $E$ is a holomorphic $r$-bundle and
$b_1(X)$ is odd,  for some line bundle $L$ with $c_1(L)=0$ and
$\deg(L,\omega)=\deg(E,\omega)/r$ the bundle
$(L^*\otimes \det E)\oplus {\cal
O}^{r-1}$ is quasi-stable, not stable, has rank and first Chern
class equal to that of $E$ and has charge $C(E)-c_2(E)$.

\medskip
\noindent{\bf Proof of Theorem~6.4:} \quad {\bf 1.}
\quad If $\{(\t X_i,\t E_i)\} \subset \b{\cal M}\backslash {\cal
M}$ is converging to $(\t X,\t E) \in {\cal M}$, there is an
irreducible component $L$ of the exceptional divisor in $\t X_i$
to which $\t E_{i_j}$ restricts non-trivially for every $i_j$ in
some subsequence. Then $H^0(L,\t E_{i_j}(-1))\not=0$ so by
semi-continuity of cohomology it follows $H^0(L,(\t E(-1))
\not=0$, implying $\t E \in \b{\cal M}\backslash {\cal M}$ and
hence that ${\cal M}$ is an open set.

If $deg(K_X,\omega)<0$ then as seen at the end of \S4, each $(\t
X,\t E)\in \b{\cal M}$ has arbitrarily small deformations
restricting trivially to the exceptional divisor in $\t X$.

\smallskip
\noindent{\bf 2,3.} \quad  Suppose $\{[(\t X_i,\t
E_i)]\} \subset \b{\cal M}$.  By Proposition~6.2, it can be
assumed that for each $i$,  $c_2(\t X_i)\le c_2(X)+N$ and by
passing to a subsequence, it can then be assumed that $c_2(\t
X_i)=c_2(X)+n$ is constant.  It is easy to see (by induction on
$n$ for example), that there is a subsequence such that $\t
X_{i_j}$ converges to a blowup $\t X$ of $X$.

There are finitely many classes in $c \in H^2(X,\m Z)^{\perp}
\subset H^2(\t X,\m Z)$ satisfying $-c\cdot c \le r^3C(E_{\rm
top})$ and therefore only finitely different classes of metrics
on $\t X_i$ of the form $\omega_{i,\alpha} :=
\pi_i^*\omega+\rho_{\alpha}$ with respect to which $\t E_i$ can
be quasi-stable for all $\alpha$ in an open set.  By passing to a
subsequence, one such class of metrics can be fixed, and then for
some new subsequence a corresponding sequence of metrics can be
found converging to a non-degenerate metric $\omega_{\alpha}$ on
$\t X$.

The $\omega_{i,\alpha}$-Hermitian-Einstein connection on $\t E_i$
has curvature bounded in $L^2(\omega_{i,\alpha})$, and hence in
$L^2(\omega_{\alpha})$ since the metrics are converging and hence
compare uniformly.  Regarding the sequence of connections as a
sequence of unitary connections on the fixed bundle $\pi^*E_{\rm
top}$ over the smooth Riemannian manifold $(X \# n\,\b{\m
P}_2,\omega_{\alpha})$ with uniformly $L^2$-bounded curvature,
Uhlenbeck's weak compactness result implies that after passing to
another subsequence there is a finite set of points where the
curvature is bubbling, with the sequence converging (after gauge
transformations) weakly in $L^2_{1,{\rm loc}}$ on the complement.
Convergence of the sequence of metrics $\{\omega_{i,\alpha}\}$
and ellipticity of the Hermitian-Einstein equations imply that
for some subsequence, this convergence can be bootstrapped to
uniform convergence in $C^2$ on compact subsets of the complement
of the bad set, converging weakly to an
$\omega_{\alpha}$-quasi-stable limit $\t E$ on $\t X$.

By the same arguments as in the proof of Theorem~1.3 of [B4],
after passing to another subsequence if necessary, there is a
further sequence of blowups $\{\t X_i'\}$ such that $\t X_i'$ is
a blowup of $\t X_i$ consisting of at most $2(C(E_{\rm top})-C(\t
E))-1$ blowups, $\pi_i'^*\t E_i$ is
$\omega_{i,\alpha,\beta}$-quasi-stable and  some sequence of
smooth integrable connections inducing $\pi'^*_i\t E_i$ converges
strongly on $\t X'$ to a bundle $\t E'$ satisfying $(\pi'_*\t
E')^{**}=\t E$.

Under the hypotheses of the third statement of the theorem,
$\pi_*\t E$ must be stable and therefore so too is $\t E$ with
respect to $\omega_{\alpha}$ for any suitable $\alpha$ with
$|\alpha|$ sufficiently small by Proposition~3.4.  Hence so too
is $\t E'$ if $|\beta|$ is small enough, proving that $\b{\cal
M}$ is sequentially compact (and hence compact, by second
countability of the topology) in this case.  Furthermore, if a
sequence $\{(\t X_{1,i},\t E_{1,i})\}\subset \b{\cal M}$
converges to $(\t X_1,\t E_1)$ and $(\t X_{1,i},\t E_{1,i})\sim
(\t X_{2,i},\t E_{2,i})$ with $\{(\t X_{2,i},\t E_{2,i})\}$
converging to $(\t X_2,\t E_2)$,   the sequence of joint blowups
$\t X_{12,i}$ converges to a joint blowup $\t X_{12}$ of $X$ and
$\pi_1^*\t E_2$, $\pi_2^*\t E_1$ must both be stable with respect
to any metric on $\t X_{12}$ of the form  $\omega_{\alpha,\beta}$
for $|\alpha|+|\beta|$ sufficiently small.  Since semi-continuity
of cohomology gives a non-zero holomorphic map between these two
pull-backs, this map must be an isomorphism.  Thus $(\t X_1,\t
E_1) \sim (\t X_2,\t E_2)$ and therefore the topology on $\b{\cal
M}$ is Hausdorff.  If $\deg(K_X,\omega)<0$ then $\deg(End_0\,\t
E\otimes K_{\t X},\omega_{\alpha})<0$ for any bundle $\t E$ on a
blowup $\t X$ stable with respect to $\omega_{\alpha}$ for
$|\alpha|$ small enough, so $H^0(\t X,End_0\,\t E) = 0 =H^2(\t
X,End_0\,\t E)$ for any such bundle, implying that the space of
deformations of $\t E$ is smooth near $\t E$.

If $b_1(X)$ is even and $rank(E_{\rm top})=2$, the same arguments
as given in \S7 of [B4] to prove Theorem~1.4 of that reference
(i.e., the last statement of Theorem~6.1 above) can be repeated
and yield the same conclusion as that above, i.e., that $\b{\cal
M}$ is compact in this case. \quad \qed

\bigskip
\bigskip
\bigskip

\centerline{\bf REFERENCES}
\bigskip
\medskip
\tolerance=1000
\parindent=25pt \everypar{\hangindent .55in}
\parskip=5pt
\frenchspacing \baselineskip=12pt

\item{[BPV]~~} W. Barth, C. Peters and A. Van de Ven, {\it
Compact Complex Surfaces}, (Springer, Berlin, Heidelberg, New
York, 1984).
\item{[BS]~~} C. B\v anic\v a  and  O. St\v an\v
asil\v a., {\it Algebraic methods in the global theory of complex
spaces}, (Wiley, London, New York, 1976).
\item{[BH]~~} P. J.
Braam and J. Hurtubise, ``Instantons on Hopf surfaces and
monopoles on solid tori", {\it J. reine angew. Math.} 400 (1989)
146--172.
\item{[Br]~~} R. Brussee, ``Stable bundles on blown up
surfaces", {\it Math. Z.} 205 (1990) 551--565.
\item{[B1]~~} N.
P. Buchdahl, ``Stable 2-bundles on Hirzebruch surfaces", {\it
Math. Z.} 194 (1987) 143--152.
\item{[B2]~~} N. P. Buchdahl,
``Hermitian-Einstein connections and stable vector bundles over
compact complex surfaces", {\it Math. Ann.} 280 (1988) 625--648.
\item{[B3]~~} N. P. Buchdahl, ``Instantons on $n\m{CP}_2$", {\it
J. Differ. Geom} 37 (1993) 669--687.
\item{[B4]~~} N. P.
Buchdahl, ``Sequences of stable vector bundles over compact
complex surfaces".  Preprint (1995).
\item{[B5]~~} N. P.
Buchdahl, ``Vector bundles on the blown-up plane".  In
preparation.
\item{[D1]~~} S. K. Donaldson, ``A new proof of a
theorem of Narasimhan and Seshadri", {\it J. Differ. Geom.} 18
(1983) 269--277.
\item{[D2]~~} S. K. Donaldson, ``Anti-self-dual
Yang-Mills connections over complex algebraic varieties and
stable vector bundles", {\it Proc. Lond. Math. Soc.} 50 (1985)
1--26.
\item{[D3]~~} S. K. Donaldson, ``Connections, cohomology
and the intersection forms of 4-manifolds", {\it J. Differ.
Geom.} 24 (1986) 275--341.
\item{[D4]~~} S. K. Donaldson,
``Irrationality and the h-cobordism conjecture ", {\it J. Differ.
Geom.} 26 (1987) 141--168.
\item{[D5]~~} S. K. Donaldson,
``Polynomial invariants for smooth four manifolds", {\it
Topology} 29 (1990) 257--315.
\item{[FS]~~} R. Fintushel and R.
J. Stern, ``The blowup formula for Donaldson invariants".
Preprint (1994).
\item{[FM1]~~} R. Friedman and J. W. Morgan,
``On the diffeomorphism types of certain algebraic surfaces I",
{\it J. Differ. Geom.} 27 (1988) 297--369.
\item{[FM2]~~} R.
Friedman and J. W.  Morgan,  ``On the diffeomorphism types of
certain algebraic surfaces II", {\it J. Differ. Geom.} 27 (1988)
371--398.
\item{[FM3]~~} R. Friedman and J. W.  Morgan, {\it
Smooth Four-Manifolds and Complex Surfaces.} Erg. Math. u.
Grenzgebiete 3 Folge 37.  Springer, Berlin Heidelberg New York,
1994.
\item{[Gau]~~} P. Gauduchon, ``Le th\'eor\`eme de
l'excentricit\'e nulle", {\it C. R. Acad. Sci. Paris} 285 (1977)
387--390.
\item{[Gie]~~} D. Gieseker, ``On the moduli of vector
bundles on an algebraic surface", {\it Ann. Math.} 106 (1977)
45--60.
\item{[K1]~~} D. Kotschick, ``On manifolds homeomorphic
to $\m{CP}^2\#8\b{\m{CP}}^2$", {\it Invent. Math.} 95 (1989)
591--600.
\item{[K2]~~} D. Kotschick, ``$SO(3)$ invariants for
$4$-manifolds with $b^+_2=1$, {\it Proc. Lond. Math. Soc.} 63
(1991) 426--448.
\item{[L]~~} M. L\"ubke, ``Chernklassen von
Hermite-Einstein Vektor-B\"undeln", {\it Math.  Ann.} 260 (1982)
133-141.
\item{[M]~~} M. Maruyama, ``On a compactification of a
moduli space of stable bundles on a rational surface". In {\it
Algebraic Geometry and Commutative Algebra}, Kinokuniya, Tokyo
1988. 233--260.
\item{[OSS]~~} C. Okonek, M. Schneider and H.
Spindler, {\it Vector bundles on complex projective spaces},
(Birkh\"auser, Boston, Basel, Stuttgart, 1980).
\item{[OV]~~} C.
Okonek and A. Van de Ven, ``Stable bundles and differentiable
structures on certain algebraic surfaces", {\it Invent. Math.} 86
(1986) 357--370.
\item{[Sch]~~} R. L. E. Schwarzenberger,
``Vector bundles on algebraic surfaces", {\it Proc. Lond. Math.
Soc. (3)} 11 (1961) 601--622.
\item{[Sed]~~} S. Sedlacek, ``A
direct method for minimizing the Yang-Mills functional", {\it
Commun. Math.  Phys.} 86 (1982) 515--528.
\item{[Ser]~~} J. -P.
Serre, ``Sur les modules projectifs.", {\it S\'em.
Dubreil--Pisot} Expos\'e 2 (1960/61) .
\item{[T]~~} C. H. Taubes,
``Self-dual Yang-Mills connections on non-self-dual 4-manifolds",
{\it J. Differ. Geom.} 17 (1982) 139--170.
\item{[U1]~~} K. K.
Uhlenbeck, ``Connections with $L^p$ bounds on curvature", {\it
Commun. Math. Phys.} 83 (1982) 31--42.
\item{[U2]~~} K. K.
Uhlenbeck, ``Removable singularities in Yang-Mills fields", {\it
Commun. Math. Phys.} 83 (1982) 11--30.

\bye